\ifdefined\PreambleLoaded\else
\documentclass{amsart}
\usepackage{marktext}
\usepackage{gimac}

\DeclareFontFamily{U}{matha}{\hyphenchar\font45}
\DeclareFontShape{U}{matha}{m}{n}{
      <5> <6> <7> <8> <9> <10> gen * matha
      <10.95> matha10 <12> <14.4> <17.28> <20.74> <24.88> matha12
      }{}
\DeclareSymbolFont{matha}{U}{matha}{m}{n}
\DeclareFontSubstitution{U}{matha}{m}{n}

\DeclareFontFamily{U}{mathx}{\hyphenchar\font45}
\DeclareFontShape{U}{mathx}{m}{n}{
      <5> <6> <7> <8> <9> <10>
      <10.95> <12> <14.4> <17.28> <20.74> <24.88>
      mathx10
      }{}
\DeclareSymbolFont{mathx}{U}{mathx}{m}{n}
\DeclareFontSubstitution{U}{mathx}{m}{n}

\DeclareMathDelimiter{\vvvert}     {0}{matha}{"7E}{mathx}{"17}
\usepackage{esint}

\newcommand{\GI}[2][]{\sidenote[colback=yellow!20]{\textbf{GI\xspace #1:} #2}}

\newcommand{\Rd}{\mathbb{R}^d}
\newcommand{\e}{\varepsilon}
\DeclareMathOperator{\Var}{Var}
\renewcommand{\var}{\Var}

\makeatletter
\newcommand\avsuminner[2]{%
	{\sbox0{$\m@th#1\sum$}%
		\vphantom{\usebox0}%
		\ooalign{%
			\hidewidth
			\smash{\,\rule[.23em]{8.8pt}{.6pt} \relax}%
			\hidewidth\cr
			$\m@th#1\sum$\cr
		}%
	}%
}
\makeatother

\makeatletter
\newcommand\avsuminnerr[2]{%
	{\sbox0{$\m@th#1\sum$}%
		\vphantom{\usebox0}%
		\ooalign{%
			\hidewidth
			\smash{\,\rule[.23em]{6pt}{0.7pt} \relax}%
			\hidewidth\cr
			$\m@th#1\sum$\cr
		}%
	}%
}
\makeatother

\newcommand{\PreambleLoaded}{}

\fi
\begin{document}
  \title[Convergence of Empirical Measures]{Convergence of Empirical Measures for i.i.d.\ samples in~$W^{-\alpha, p}$}
  \begin{abstract}
    Given~$N$ i.i.d.\ samples from a probability measure~$\mu$ on~$\R^d$, we study the rate of convergence of the empirical measure~$\mu_N \to \mu$ in the negative Sobolev space~$W^{-\alpha, p}$.
    When~$W^{-\alpha, p}$ contains point measures (i.e.\ when~$\alpha p > (p-1)d$), we show~$\E \norm{\mu_N - \mu}_{W^{-\alpha, p}}^p \leq C_d / N^{p/2}$ for an explicit dimensional constant~$C_d$, and obtain a Gaussian tail bound.
    When~$0 < \alpha p \leq d(p-1)$, we prove a similar result for Gaussian regularizations.
  \end{abstract}
  \author[Iyer]{Gautam Iyer}
  \address{%
    Department of Mathematical Sciences, Carnegie Mellon University, Pittsburgh, PA 15213.
  }
  \email{gautam@math.cmu.edu}
  \author[Venkatraman]{Raghavendra Venkatraman}
  \address{%
    Department of Mathematics, University of Utah, Salt Lake City, UT 84112.
  }
  \email{raghav@math.utah.edu}
  \thanks{This work has been partially supported by the National Science Foundation under grants
    DMS-2342349, 
    DMS-2406853, 
    DMS-2407592. 
    The authors also thank the SIAM Northern States Section and the Center for Nonlinear Analysis.}
  
  \maketitle
\section{Introduction.}\label{s:intro}

\subsection{Main Results}
\label{ss.main}
Let~$\mu$ be a probability measure on~$\R^d$, and $X_1$, \dots, $X_N$ be~$N$ i.i.d.\ samples from~$\mu$.
Let~$\mu_N$ be the empirical measure
\begin{equation}
  \mu_N \defeq \frac{1}{N} \sum_{i=1}^N \delta_{X_i}
  \,,
\end{equation}
where~$\delta_x$ denotes the Dirac~$\delta$-measure at~$x$.
The aim of this paper is to study the convergence of~$\mu_N \to \mu$ in negative Sobolev spaces.
For~$L^2$ based Sobolev spaces, $\E \norm{\mu_N - \mu}_{H^{-s}}^2$ can be computed exactly (see Proposition~\ref{p:EmpMeasureHAlpha}, below).
The main result of this paper bounds the convergence rate of~$\mu_N \to \mu$ in~$L^p$ based Sobolev spaces for all~$p \in (1, \infty)$.

To state our result, let~$\alpha > 0$, $p \in (1, \infty)$ and~$ W^{-\alpha,p} =W^{-\alpha, p}(\Rd) $ be the~$L^p$ based Sobolev space with index~$-\alpha$.
While there are many equivalent definitions of a norm in~$W^{-\alpha, p}$, the version that is most convenient for us is based on the heat kernel, and stated precisely in Section~\ref{s:main}, below.
We show that when~$\alpha p > (p-1)d$ (i.e.\ when the space~$W^{-\alpha, p}$ contains point measures), we have $\E \norm{\mu_N - \mu}_{W^{-\alpha, p}}^p \leq O(1/N^{p/2})$ and a Gaussian tail bound.

\begin{theorem}\label{t:mainIntro1}
  Let~$p \in (1, \infty)$, $q = p / (p-1)$ be its H\"older conjugate, and suppose~$\alpha > d/q = \alpha p / (p-1)$.
  There exists an absolute constant~$C$ (independent of~$\mu$, $p$ and~$d$) such that for every~$N \in \mathbb{N}$ we have
  \begin{equation}\label{e:muNMuW}
    \paren[\big]{\E \norm{\mu_N - \mu}_{W^{-\alpha, p}}^p}^{1/p}
      \leq
	\frac{C d^{1/p}}{\sqrt{N}} \norm{\delta_0}_{W^{-\alpha, p}}
	\,.
  \end{equation} 
  Here~$\delta_0$ is the Dirac~$\delta$-distribution at~$0$.
  Additionally, for every~$N \in \N$, the random variable~$\norm{\mu_N - \mu}_{W^{-\alpha, p}}$ satisfies the Gaussian tail bound
  \begin{equation}
    \P\paren[\big]{
      \norm{\mu_N - \mu}_{W^{-\alpha, p}} > \lambda
    }
    \leq 2 \exp\paren[\bigg]{ \frac{-N \lambda^2}{C d^{2/p} p \norm{\delta_0}_{W^{-\alpha, p}}^2 } }
    \,,
  \end{equation}
  for every~$\lambda > 0$.
\end{theorem}

Notice that when~$\alpha \leq d/q$, the empirical measure~$\mu_N \not\in W^{-\alpha, p}$, and so both the left and right hand side of~\eqref{e:muNMuW} are infinite.
In this case it is natural to study the convergence of the Gaussian regularized point measures instead.
For any~$\epsilon > 0$, define~$\mu^\epsilon \defeq \mu * \varPhi_\epsilon$ and $\mu_N^\epsilon \defeq \mu_N * \varPhi_\epsilon$, where~$\varPhi_\epsilon$ is the density of a centered radial Gaussian with variance~$2\epsilon$ (see~\eqref{e.heatkernel}, below).
When~$\alpha \leq d/q$ we bound~$\E \norm{\mu^\epsilon_N - \mu}_{W^{-\alpha, p}}^p$, and obtain a Gaussian tail bound.
\begin{theorem}\label{t:mainRegularized}
  Let~$p \in (1, \infty)$, $\alpha > 0$ and suppose~$\mu \in W^{-\alpha, p}$.
  There exists an absolute constant~$C$ (independent of~$\mu$, $p$ and~$d$) such that for every~$N \in \mathbb{N}$, we have
  \begin{equation}\label{e:MuNEpMuW}
    \paren[\big]{\E \norm{\mu^\epsilon_N - \mu}_{W^{-\alpha, p}}^p}^{1/p}
      \leq
	\frac{C d^{1/p}}{\sqrt{N}} \norm{\varPhi_\epsilon}_{W^{-\alpha, p}}
	+ \norm{\mu^\epsilon - \mu}_{W^{-\alpha, p}}
	\,.
  \end{equation} 
  Additionally, for every~$N \in \N$, $\lambda > 0$ we have the Gaussian tail bound
  \begin{equation}
    \label{e:MuNEpMuConcentration}
    \P\paren[\big]{
      \norm{\mu^\epsilon_N - \mu}_{W^{-\alpha, p}} > \lambda
    }
    \leq 2 \exp\paren[\bigg]{ \frac{-N \paren[\big]{\lambda - \norm{\mu^\epsilon - \mu}_{W^{-\alpha, p}} }_{+}^2}{C p d^{2/p} \norm{\varPhi_\epsilon}_{W^{-\alpha, p}}^2 } }
    \,.
  \end{equation}
\end{theorem}

In order to use Theorem~\ref{t:mainRegularized} one needs to understand the behavior of~$\norm{\varPhi_\epsilon}_{W^{-\alpha, p}}$ and~$\norm{\mu^\epsilon - \mu}_{W^{-\alpha, p}}$ as~$\epsilon \to 0$.
These are both well known and stated here as remarks for easy reference.
\GI[2025-12-10]{Find a reference for both.}

\begin{remark}
  For every~$\mu \in W^{-\alpha, p}$, we know~$\norm{\mu^\epsilon - \mu}_{W^{-\alpha, p}}$ vanishes with~$\epsilon$, but the rate may be arbitrarily slow in~$\epsilon$.
  If, however, $\mu$ was more regular than~$W^{-\alpha, p}$, then one has an explicit rate.
  For instance, if $\mu \in W^{-\alpha', p}$ for some~$\alpha' < \alpha$, then standard regularity estimates imply
  \begin{equation}
    \norm{\mu^\epsilon - \mu}_{W^{-\alpha, p}} \leq C_{\alpha, d} \epsilon^{\frac{\alpha - \alpha'}{2} \varmin 1}
    \,.
  \end{equation}
\end{remark}

\begin{remark}\label{r:Gaussian}
  It is well known that
  \begin{equation}
    \norm{\varPhi_\epsilon}_{W^{-\alpha, p}} \approx
      \begin{dcases}
	 \abs{\ln \epsilon}^{1/p}
	  & \alpha = \tfrac{d}{q}
	\,
      \\
	\frac{1}{ \epsilon^{\frac{d}{q} - \alpha}}
	  & \alpha < \tfrac{d}{q}
	\,.
      \end{dcases}
  \end{equation}
  with constants that depend on~$\alpha$, $p$, and~$d$.
  For convenience, we state this precisely in the appendix (see Proposition~\ref{p:PhiEpBounds}), and present a proof with explicit bounds on the constants.
\end{remark}

The constant~$d^{1/p}$ appearing in~\eqref{e:muNMuW} can be replaced with the~$p$-th root of the best weak~$L^1$ bound for the Hardy--Littlewood maximal function.
The standard proof using the Vitali covering lemma shows this constant is at most~$3^d$, and a classical result of Stein and Str\"omberg~\cite{SteinStromberg83} improves this to~$O(d)$.
To the best of our knowledge the sharp constant for this weak~$L^1$ bound is not known.

Finally we mention that our results do not require any regularity or moment conditions on the probability measure~$\mu$,  nor any bounds from below for~$\mu$ on its support.
Moreover a standard localization argument can be used to obtain similar rates of convergence when~$\mathbb{R}^d$ is replaced by a closed manifold.
We also note that without further assumptions matching lower bounds for~\eqref{e:muNMuW} can not be obtained, as can be seen by taking~$\mu=\delta_0$.

\subsection{Motivation and literature review}\label{s:lit}

The convergence of the empirical measure~$\mu_N \to \mu$ is a central question in probability and statistics, and has been extensively studied.
The law of large numbers shows that almost surely, the sequence~$\mu_N \to \mu$ weakly in the space of probability measures.

There are a number of results quantifying the weak convergence of~$\mu_N \to \mu$.
One popular approach is to use the \emph{Wasserstein distance}, which is natural as it measures the transportation cost of moving~$\mu$ to~$\nu$.
Explicitly for~$p \in [1, \infty)$, the \emph{Wasserstein distance} between two probability measures~$\mu, \nu$ (see for instance~\cite{Villani09a}) is defined by
\begin{equation}
  W_p(\mu, \nu)
    = \inf_{\gamma \in \Gamma(\mu, \nu)} \paren[\big]{\E_{(X, Y) \sim \gamma} \abs{X - Y}^p}^{1/p}
    \,,
\end{equation}
where~$\Gamma(\mu, \nu)$ is the set of all couplings of~$\mu$ and~$\nu$. 

Seminal work of Ajtai, Kolm\'os, Tusn\'ady~\cite{AjtaiKomlosEA84} studies the rate of convergence of~$W_p(\mu_N, \mu)$ in the context of an optimal matching problem when~$\mu = \operatorname{Unif}([0,1]^d)$ is the uniform distribution.
It has since been studied by several authors~\cite{Talagrand21,Ledoux23}, and recent work of Goldman and Trevisan~\cite{GoldmanTrevisan21} shows
\begin{equation}
  \lim_{N \to \infty} N^{1/d} \paren[\big]{\E W_p^p(\mu_N, \mu)}^{1/p} \in (0, \infty)
  \,,
  \quad\text{when }
  \mu = \operatorname{Unif}([0,1]^d)\,,~
  d \geq 3
  \,.
\end{equation}
This gives both an upper and lower bound showing $(\E W_p^p(\mu_N, \mu))^{1/p} \approx 1/N^{p/d}$.

For general distributions, to the best of our knowledge, the sharp rate of convergence is not known.
Under a suitable moment condition, one can show~\cite{
  DereichScheutzowEA13,
  BoissardLeGouic14,
  FournierGuillin15
} the upper bound
\begin{equation}
  (\E W^p_p(\mu_N, \mu))^{1/p} \leq \frac{C}{N^{1/d}}
    \quad\text{when } 1 \leq p < \frac{d}{2}\,.
\end{equation}
This has been extended in~\cite{LedouxZhu21} for Gaussians for all~$p \in [1, d)$, and by Divol~\cite{Divol21} for all~$p \in [1, \infty)$ for compactly supported distributions with strictly positive densities (see~\cite{LeightonShor89,GarciaTrillosSlepcev15} for the~$p = \infty$ case).
More generally  Weed and Bach~\cite{WeedBach19} show the bounds
\begin{equation}
  \E W_p(\mu_N, \mu) \leq \frac{C}{N^{1/t}}
  \quad\text{and}\quad
  W_p(\mu_N, \mu) \geq \frac{C}{N^{1/s}}
\end{equation}
for every~$t > d^*_p(\mu)$ and~$s < d_*(\mu)$.
Here~$d^*_p(\mu)$ and~$d_*(\mu)$ are the upper and lower Wasserstein dimensions respectively (see Definition~4 in~\cite{WeedBach19}).
(Note the lower bound holds surely, and not on average.)

The Wasserstein bounds above suffer from the curse of dimensionality -- the convergence rate becomes extremely slow when the dimension is large.
Our aim in this paper is to look for results that allow for faster convergence in high dimensions.
To hint at why this is to be expected, we note that for any (nice enough) test function~$f$, we have the \emph{dimension independent} identity
\begin{equation}
  \E \paren[\bigg]{ \int_{\R^d} f \, d(\mu_N - \mu)  }^2
    = \frac{\var_\mu(f)}{N}
    \,,
\end{equation}
which is the foundation of several Monte Carlo methods used today (see for instance~\cite{Liu08}).
Of course, for the above to be defined, the function~$f$ needs to be regular enough to allow for evaluation at points.
This suggests dimension independent bounds may require the use of spaces which discount fine scale structure.

One such example is a maximum mean discrepancy (MMD) distance, which is known to quantify weak convergence~\cite{SimonGabrielBarpEA23}.
With this approach, the probability measures are embedded into a Hilbert space and the norm in this space is used to measure distance.
In this setting the expected value of the MMD distance squared can be computed exactly using the modulated energy~\cite{HessChildsRosenzweigEA25} (see also Proposition~\ref{p:EmpMeasureHAlpha}, below).

The present paper aims to study convergence rates outside the Hilbert space setting, where an exact calculation is not available.
The simplest and most natural example is to consider \emph{negative Sobolev spaces}, which have been used by many authors to quantify weak convergence~\cite{Thiffeault12}.
Here we obtain rates of convergence with a dimension independent exponent.
More precisely, Theorem~\ref{t:mainIntro1} shows the convergence rate in~$L^p$ based Sobolev spacesis~$O(C_d / \sqrt{N})$, for an explicit dimensional constant~$C_d$.
The disadvantage of using negative Sobolev spaces is that for some indices the space may not contain point measure~$\mu_N$, or even the target measure~$\mu$.
We address this by studying the convergence of Gaussian regularized point measures instead, as stated in Theorem~\ref{t:mainRegularized}.


\subsection*{Plan of this paper}

In Section~\ref{s:main} we state the main result of this paper (Theorem~\ref{t:empMeas}) which yields both Theorem~\ref{t:mainIntro1} and~\ref{t:mainRegularized}.
We prove Theorem~\ref{t:empMeas} in Section~\ref{s:calW} when~$\alpha \not\in \N$, and in Section~\ref{s:scrW} when~$\alpha \in \N$.
In Section~\ref{s:Halpha} we prove Proposition~\ref{p:EmpMeasureHAlpha}, which is a stronger version of Theorem~\ref{t:empMeas} for~$p = 2$.
Finally in the appendix we compute the negative Sobolev norms of Gaussians, as stated in Remark~\ref{r:Gaussian}

\subsection*{Acknowledgements}

The authors thank
  Giovanni Leoni,
  Dejan Slep\v cev,
  and Matthew Rosenzweig
for helpful discussions.

\section{Main results and proofs.}\label{s:main}

Before proving Theorems~\ref{t:mainIntro1} and~\ref{t:mainRegularized}, we define the~$W^{-\alpha, p}$ norm that is used in the statement of the theorem.
There are several equivalent definitions of the norm in~$W^{-\alpha, p}$ (see for instance~\cite{Triebel92}), but the version that is most convenient for us is the one from Appendix~D in~\cite{ArmstrongKuusiEA19} based on heat kernels.
For~$\alpha > 0$, $p \in (1, \infty)$ we define the~$W^{-\alpha, p}$ norm by
\begin{equation}\label{e:normW}
  \norm{u}_{W^{-\alpha, p}} =
    \begin{dcases}
      \norm{u}_{\mathcal W^{-\alpha, p}} & \alpha \not\in \N\,,
      \\
      \norm{u}_{\mathscr W^{-\alpha, p}} & \alpha \in \N\,.
    \end{dcases}
\end{equation}
Here the norms~$\norm{\cdot}_{\mathcal W^{-\alpha, p}}$ and~$\norm{\cdot}_{\mathscr W^{-\alpha, p}}$ are defined by
\begin{subequations}
\begin{align}
  \label{e:WMinusAlphaNotInt}
  \norm{u}_{\mathcal W^{-\alpha, p}}
    &\defeq \paren[\Big]{ \int_0^1 \paren[\big]{ t^{\frac{\alpha}2} \|u * \varPhi_t\|_{L^p(\Rd)}}^p \frac{\,dt}{t} }^{\frac1{p}}
    \,,
    \\
    \label{e:WMinusAlphaInt}
  \norm{u}_{\mathscr W^{-\alpha, p}}
  &\defeq
  \paren[\bigg]{ \int_{\Rd} \paren[\Big]{ \int_0^1 \paren[\big]{ t^{\frac{\alpha}2} |u * \varPhi_t|}^2  \frac{\,dt}{t}}^{\frac{p}2}  \,dx }^{\frac1{p}}
  \,,
\end{align}
\end{subequations}
and~$\varPhi$ is the standard Euclidean heat kernel
\begin{equation}\label{e.heatkernel}
  \varPhi_t (x) \defeq \frac{1}{(4\pi t)^{d/2}} \exp \paren[\Big]{ - \frac{|x|^2}{4t}}\,,
  \quad x \in \mathbb{R}^d\,,~t >0\,,
\end{equation}
For convenience, when~$t = 0$ we will use~$\varPhi_0$ to denote the Dirac~$\delta$ distribution at~$0$.
Appendix~D in~\cite{ArmstrongKuusiEA19} provides a detailed proof showing~\eqref{e:normW} is equivalent to the standard~$W^{-\alpha, p}$ norm based on duality.

Both Theorems~\ref{t:mainIntro1} and~\ref{t:mainRegularized} follow immediately from the following result for Gaussian regularizations.

\begin{theorem}\label{t:empMeas}
  There exists an absolute constant~$C$, independent of~$d$ and~$\mu$, such that for any~$p\in (1,\infty)$ and~$\epsilon \geq 0$ we have the moment bound
  \begin{equation}\label{e:MuNMu}
    \paren[\big]{\E \norm{\mu^\epsilon_N - \mu^\epsilon}_{W^{-\alpha, p}}^p}^{1/p}
      \leq
	\frac{C d^{1/p} \sqrt{p}}{\sqrt{N}} \norm{\varPhi_\epsilon}_{W^{-\alpha, p}}
	\,,
  \end{equation}
  and the Gaussian tail bound
  \begin{equation}\label{e:MuNMuSubGW}
    \P\paren[\big]{
      \norm{\mu^\epsilon_N - \mu^\epsilon}_{W^{-\alpha, p}} > \lambda
    }
    \leq 2 \exp\paren[\bigg]{ \frac{-N \lambda^2}{C d^{2/p} p \norm{\varPhi_\epsilon}_{W^{-\alpha, p}}^2 } }
    \,,
    \quad\text{for every } \lambda > 0\,.
  \end{equation}
  If~$\alpha p > (p - 1) d$ and~$\mu \in W^{-\alpha, p}$, then both~\eqref{e:MuNMu} and~\eqref{e:MuNMuSubGW} also hold for~$\epsilon = 0$.
\end{theorem}

Theorems~\ref{t:mainIntro1} and~\ref{t:mainRegularized} follow immediately from Theorem~\ref{t:empMeas}.
\begin{proof}[Proof of Theorems~\ref{t:mainIntro1} and~\ref{t:mainRegularized}]
  When~$\alpha p > d/ q$, Theorem~\ref{t:empMeas} with~$\epsilon = 0$ reduces to Theorem~\ref{t:mainIntro1}.
  For Theorem~\ref{t:mainRegularized}, the triangle inequality implies
  \begin{equation}
    \paren[\big]{\E \norm{\mu^\epsilon_N - \mu}_{W^{-\alpha, p}}^p}^{1/p}
      \leq
      \paren[\big]{\E \norm{\mu^\epsilon_N - \mu^\epsilon}_{W^{-\alpha, p}}^p}^{1/p}
      + \norm{\mu^\epsilon - \mu}_{W^{-\alpha, p}}
      \,.
  \end{equation} 
  Now~\eqref{e:MuNEpMuW} and~\eqref{e:MuNEpMuConcentration} follow from~\eqref{e:MuNMu} and~\eqref{e:MuNMuSubGW} respectively.
\end{proof}

We prove Theorem~\ref{t:empMeas} in Section~\ref{s:calW} when~$\alpha \not \in \N$, and in Section~\ref{s:scrW} when~$\alpha \in \N$, and these are the main contributions of this paper.
To briefly comment on the proof, we note that the standard dual norm on~$W^{-\alpha, p}$ is not convenient to work with in this context, which is why we use the heat kernel based bounds~\eqref{e:WMinusAlphaNotInt}--\eqref{e:WMinusAlphaInt}.
As we will see shortly, the main tool used in the proof is a bound on the \emph{subgaussian norm} of the random variable~$\varPhi_t(x - X_1)$, treating~$t \in (0, 1)$ and~$x \in \R^d$ as parameters.
Of course, since~$\varPhi_t \in L^\infty$, the subgaussian norm is easily bounded by~$\norm{\varPhi_t}_{L^\infty}$.
However, this bound is too crude to yield Theorem~\ref{t:empMeas}.
What we need is a bound that decays nicely in \emph{both} the parameters~$x$ and~$t$.
We prove this bound in Lemma~\ref{l:subGbound} (below) and control the subgaussian norm in terms of the Hardy--Littlewood maximal function of~$\mu$.
We then use the weak-$L^1$ bound on the maximal function to prove Theorem~\ref{t:empMeas}.

Before delving into the proof, we briefly digress and state the analog of Theorem~\ref{t:empMeas} when~$p = 2$.
In this case the expressions in both~\eqref{e:WMinusAlphaNotInt} and~\eqref{e:WMinusAlphaInt} agree, and there is a slightly stronger result that can be obtained with a simpler proof.
To simplify notation let~$H^{-\alpha}$ denote the $L^2$-based Sobolev space with index~$-\alpha$, and define
\begin{equation}\label{e:HMinusAlpha}
  \norm{u}_{H^{-\alpha}}
    \defeq
      \norm{u}_{\mathcal W^{-\alpha, 2}}=
      \norm{u}_{\mathscr W^{-\alpha, 2}}
      = \paren[\Big]{ \int_0^1 \paren[\big]{ t^{\frac{\alpha}2} \|u * \varPhi_t\|_{L^2(\Rd)}}^2 \frac{\,dt}{t} }^{\frac1{2}}
    \,,
\end{equation}
whether or not~$\alpha \in \N$.
Our next result shows that when~$p = 2$ the bound~\eqref{e:MuNMu} can be made exact equality, and~\eqref{e:MuNMuSubGW} can be proved without any dimensionally dependent constants.

\begin{proposition}\label{p:EmpMeasureHAlpha}
  For every~$\alpha > 0$, $\epsilon > 0$ and every~$N \in \N$ we have
  \begin{equation}\label{e:MuNMinusMuEpHAlpha}
    \E \|\mu_N^\e  - \mu^\e \|_{H^{-\alpha}}^2 
      = \frac{1}N \paren[\Big]{
	\norm{\varPhi_\epsilon}_{H^{-\alpha}}^2
	  - \|\mu_\e\|_{H^{-\alpha}}^2
	} \,.
  \end{equation}
  Additionally, we have the Gaussian tail bound
  \begin{equation}\label{e:MuNMuSubGHalpha}
    \P\paren[\big]{
      \norm{\mu^\epsilon_N - \mu^\epsilon}_{H^{-\alpha}} > \lambda
    }
    \leq 2 \exp\paren[\bigg]{ \frac{-N \lambda^2}{C \norm{\varPhi_\epsilon}_{H^{-\alpha}}^2 } }
    \,,
    \quad\text{for every } \lambda > 0\,.
  \end{equation}
  If~$\alpha > \frac{d}2$ and~$\mu \in H^{-\alpha}$ then both~\eqref{e:MuNMinusMuEpHAlpha} and~\eqref{e:MuNMuSubGHalpha} also hold for~$\epsilon = 0$.
\end{proposition}


\section{Convergence in \texorpdfstring{$\mathcal W^{-\alpha, p}$}{cal W -alpha p}.}\label{s:calW}

The aim of this section is to prove Theorem~\ref{t:empMeas} when~$\alpha \not \in \N$.
We do this by proving a slightly more general result.
Let~$\mathcal W^{-\alpha , p} =\mathcal W^{-\alpha , p}(\Rd) $ be the space of all tempered distributions for which the norm in~\eqref{e:WMinusAlphaNotInt} is finite.
We will now show that Theorem~\ref{t:empMeas} holds with~$\norm{\cdot}_{W^{-\alpha, p}}$ replaced with~$\norm{\cdot}_{\mathcal W^{-\alpha, p}}$, whether or not~$\alpha \in \N$.
Explicitly, the result we prove is as follows.

\begin{proposition}\label{p:empMeasCalW}
  There exists an absolute constant~$C$, independent of~$d$ and~$\mu$, such that for any~$p\in (1,\infty)$ and~$\epsilon \geq 0$ we have the moment bound
  \begin{equation}\label{e:MuNMuCalW}
    \paren[\big]{\E \norm{\mu^\epsilon_N - \mu^\epsilon}_{\mathcal W^{-\alpha, p}}^p}^{1/p}
      \leq
	\frac{C d^{1/p} \sqrt{p}}{\sqrt{N}} \norm{\varPhi_\epsilon}_{\mathcal W^{-\alpha, p}}
	\,,
  \end{equation}
  and the Gaussian tail bound
  \begin{equation}\label{e:MuNMuSubG}
    \P\paren[\big]{
      \norm{\mu^\epsilon_N - \mu^\epsilon}_{\mathcal W^{-\alpha, p}} > \lambda
    }
    \leq 2 \exp\paren[\bigg]{ \frac{-N \lambda^2}{C d^{2/p} p \norm{\varPhi_\epsilon}_{\mathcal W^{-\alpha, p}}^2 } }
    \,,
    \quad\text{for every } \lambda > 0\,.
  \end{equation}
  If~$\alpha p > (p - 1) d$ and~$\mu \in \mathcal W^{-\alpha, p}$, then both~\eqref{e:MuNMuCalW} and~\eqref{e:MuNMuSubG} also hold for~$\epsilon = 0$.
\end{proposition}

We state this version here for convenience, and note that it immediately implies Theorem~\ref{t:empMeas} for~$\alpha \not\in \N$.

\subsection{Proof of Proposition~\ref{p:empMeasCalW}}

As mentioned earlier, the main idea behind the proof of Proposition~\ref{p:empMeasCalW} is to consider the random variable~$\varPhi_t(x - X_1)$, treating~$x \in \R^d$ and $t \in (0, 1)$ as parameters, and show that the \emph{subgaussian norm} decays suitably in both~$x$ and~$t$.
Recall (see for instance Section~2.5 in~\cite{Vershynin18}) that if~$Z$ is a one dimensional random variable its subgaussian norm, denoted by~$\norm{Z}_{\psi_2}$ is defined by
\begin{equation}
  \norm{Z}_{\psi_2} = \inf\set[\Big]{
    c > 0 \st \E \exp\paren[\Big]{\frac{Z^2}{c^2} \leq 2}
  }
  \,.
\end{equation}
What we need to prove Theorem~\ref{t:empMeas} is to show that $\norm{\varPhi_t(x - X_1)}_{\psi_2}$ essentially decays like a Gaussian in~$x, t$.
This is our first lemma.

\begin{lemma}\label{l:subGbound}
  For any~$t\in (0,1), x \in \Rd$ we have  
  \begin{equation}\label{e:subGbound}
    \sigma(x, t) \defeq \|\varPhi_t(x-X_1)\|_{\Psi_2}
      \leq \frac{C}{(4\pi t)^{d/2}} \exp \paren[\bigg]{ \frac{-1}{ 4t \paren[\big]{V_d M \mu(x)}^{2/d}} }\,,
  \end{equation}
  where and~$M\mu$ denotes the Hardy--Littlewood maximal function of~$\mu$.
  Here~$V_d$ denotes the volume of the unit ball in~$\Rd$.
\end{lemma}

Here clarify~$\norm{\varPhi_t(x - X_1)}_{\psi_2}$ denotes the subgaussian norm of the random variable~$\varPhi_t(x - X_1)$, where~$x, t$ are fixed, deterministic, parameters.
We prove Lemma~\ref{l:subGbound} in Section~\ref{s:subGProof}, below.

To see how~\eqref{e:subGbound} decays in~$x, t$, we note that the maximal function of any probability measure ``usually'' decays
like~$1/\abs{x}^d$ as~$\abs{x} \to \infty$.
When this happens,~$\sigma(x, t)$ will satisfy a Gaussian upper bound.
Of course, the maximal function does not always decay like~$1/\abs{x}^d$, and obtaining the decay we need from~\eqref{e:subGbound} requires a little work.
We do this by combining~\eqref{e:subGbound} with the weak~$L^1$ bound for the maximal function.
Explicitly, what we need is the following lemma.

\begin{lemma}\label{l:sigmaCalW}
  There exists an absolute constant~$C$, independent of~$d$, $\mu$, such that for every~$p \in (1, \infty)$, $\alpha > 0$, $\epsilon \geq 0$ we have
  \begin{equation}\label{e:intSigmaBd}
    \paren[\Big]{
      \int_0^1 
	t^{\frac{\alpha p}{2} - 1} \norm{\sigma(\cdot, t + \epsilon)}_{L^p}^p
	\, dt
    }^{1/p}
      \leq C d^{1/p} \norm{\varPhi_\epsilon}_{\mathcal W^{-\alpha, p}}
    %
      \,.
  \end{equation}
\end{lemma}
\begin{remark}\label{r:alphaBigCalW}
  Note that when~$\alpha > d/q$ (i.e.\ when~$\alpha p > (p-1) d$), the right hand side of~\eqref{e:intSigmaBd} is finite when~$\epsilon = 0$.
\end{remark}

Momentarily postponing the proofs of Lemmas~\ref{e:MuNMuSubG} and~\ref{l:sigmaCalW} we now prove Proposition~\ref{p:empMeasCalW}.
Here and throughout this paper we use $C>0$ to denote an absolute constant that does not depend on the measure~$\mu$, or the parameters~$p>1,d\in \mathbb{N},\alpha>0$, and may increase from line to line.

\begin{proof}[Proof of Proposition~\ref{p:empMeasCalW}]
  For any~$x \in \R^d$, $t \in \R$, define the random variable~$S_N(x, t)$ by
  \begin{equation}\label{e:SNdef}
    S_N(x, t)
      \defeq (\mu_N - \mu) * \varPhi_t(x)
      = \frac{1}{N} \sum_{i=1}^N \paren[\big]{ \varPhi_t(x - X_i) - \mu_t(x) }
    \,.
  \end{equation}
  Now we compute
  \begin{equation}\label{e:MuNFubini}
    \E \norm{\mu^\epsilon_N - \mu^\epsilon}_{\mathcal W^{-\alpha, p}}^p
      = \int_{t = 0}^1 \int_{x\in \R^d} t^{\frac{\alpha p}{2} - 1} \E \abs{S_N(x, t+\epsilon)}^p \, dx \, dt
      \,.
  \end{equation}
  Since the~$p$-th moment of a random variable is controlled by~$\sqrt{p}$ times its subgaussian norm, this implies
  \begin{equation}\label{e:MuN1}
    \E \norm{\mu^\epsilon_N - \mu^\epsilon}_{\mathcal W^{-\alpha, p}}^p
      \leq C^p p^{p/2} \int_{t = 0}^1 \int_{x\in \R^d} t^{\frac{\alpha p}{2} - 1} \norm{S_N(x, t+\epsilon)}_{\psi_2}^p \, dx \, dt
      \,.
  \end{equation}

  Now by Lemma~\ref{l:subGbound} and Hoeffding's inequality the random variable~$S_N(x, t+\epsilon)$ is subgaussian with norm controlled by
  \begin{equation}\label{e:SnPsi2}
    \norm{S_N(x, t+\epsilon)}_{\psi_2} \leq \frac{C \sigma(x, t+\epsilon)}{\sqrt{N}}
    \,.
  \end{equation}
  Using this in~\eqref{e:MuN1} implies
  \begin{equation}
    \paren[\big]{\E \norm{\mu^\epsilon_N - \mu^\epsilon}_{\mathcal W^{-\alpha, p}}^p}^{1/p}
      \leq \frac{C \sqrt{p}}{\sqrt{N}}
	\paren[\Big]{
	  \int_0^1 
	    t^{\frac{\alpha p}{2} - 1} \norm{\sigma(\cdot, t + \epsilon)}_{L^p}^p
	    \, dt
	}^{1/p}
      \,.
  \end{equation}
  Combining this with Lemma~\ref{l:sigmaCalW} yields~\eqref{e:MuNMuCalW} as desired.
  \medskip

  We now prove the Gaussian tail bound~\eqref{e:MuNMuSubG}.
  Notice
  \begin{equation}
    \norm{\mu^\epsilon_N - \mu^\epsilon}_{\mathcal W^{-\alpha, p}}
      = f_N(X_1, \dots, X_N)
  \end{equation}
  where
  \begin{align}
    f_N(y_1, \dots, y_N)
      &\defeq
	\norm{g_N(\cdot, \cdot; y_1, \dots, y_N)}_*
    \\
    \label{e:gNdef}
    g_N(x, t; y_1, \dots, y_N)
      &\defeq
      \frac{1}{N} \sum_{i=1}^N \varPhi_{t+\epsilon}(x - y_i) - \mu_{t+\epsilon}(x)
      \,,
  \end{align}
  and~$\norm{\cdot}_*$ is the weighted~$L^p$ norm
  \begin{equation}
    \norm{h}_*
      \defeq
	\paren[\bigg]{
	  \int_{0}^1 \int_{\R^d} t^{\frac{\alpha p}{2} - 1}
	    \abs{h(x, t)}^p \, dx \, dt
	}^{1/p}
	\,.
  \end{equation}

  We now claim for any~$i \in \set{1, \dots, N}$, $y_i, y_i' \in \R^d$ the function~$f_N$ satisfies
  \begin{equation}\label{e:bddDiff}
    \abs{f_N(y_1, \dots, y_N) - f_N(y_1, \dots, y_{i-1}, y_i', y_{i+1}, \dots, y_N)}
      \leq \frac{2}{N} \norm{\varPhi_\epsilon}_{\mathcal W^{-\alpha, p}}
      \,.
  \end{equation}
  To prove this we note that by symmetry it suffices to prove~\eqref{e:bddDiff} for~$i = 1$.
  Let~$y_1, y_1' \in \R^d$, and~$\tilde y = (y_2, \dots, y_N)$.
  Since~$\norm{\cdot}_*$ is a norm, the reverse triangle inequality implies
  \begin{equation}
    \abs{f_N(y_1, \tilde y) - f_N(y_1', \tilde y)}
      \leq \norm{h(\cdot, \cdot; y_1, y_1')}_*
      \,,
  \end{equation}
  where
  \begin{equation}\label{e:hNDef}
    h(x, t; y_1, y_1')
      \defeq \frac{1}{N} \paren[\big]{
	\varPhi_{t+\epsilon}(x - y_1) - \varPhi_{t+\epsilon}(x - y_1') 
      }
      \,.
  \end{equation}
  Clearly
  \begin{equation}
    \norm{h(\cdot, \cdot; y_1, y_1')}_{*}
      \leq \frac{2}{N} \norm{\varPhi_\epsilon}_{\mathcal W^{-\alpha, p}}
  \end{equation}
  from which~\eqref{e:bddDiff} follows.

  Since the bounded differences property~\eqref{e:bddDiff} holds, we may now apply McDiarmid's inequality to show~$\norm{\mu^\epsilon_N - \mu^\epsilon}_{\mathcal W^{-\alpha, p}}$ is subgaussian with
  \begin{equation}
    \norm[\Big]{
      \norm{ \mu^\epsilon_N - \mu^\epsilon}_{\mathcal W^{-\alpha, p}}
       - \E \norm{ \mu^\epsilon_N - \mu^\epsilon}_{\mathcal W^{-\alpha, p}}
    }_{\psi_2}
    \leq \frac{C}{\sqrt{N}} \norm{\varPhi_\epsilon}_{\mathcal W^{-\alpha, p}}
    \,.
  \end{equation}
  Thus
  \begin{align}
    \norm[\Big]{
      \norm{ \mu^\epsilon_N - \mu^\epsilon}_{\mathcal W^{-\alpha, p}}
    }_{\psi_2}
    &\leq 
      \frac{C}{\sqrt{N}} \norm{\varPhi_\epsilon}_{\mathcal W^{-\alpha, p}}
       + C \E \norm{ \mu^\epsilon_N - \mu^\epsilon}_{\mathcal W^{-\alpha, p}}
  \\
    &\leq 
      \frac{C}{\sqrt{N}} \norm{\varPhi_\epsilon}_{\mathcal W^{-\alpha, p}}
       + C \paren[\big]{\E \norm{ \mu^\epsilon_N - \mu^\epsilon}_{\mathcal W^{-\alpha, p}}^p}^{1/p}
    \,.
  \end{align}
  Using~\eqref{e:MuNMuCalW} this implies~$\norm{ \mu^\epsilon_N - \mu^\epsilon}_{\mathcal W^{-\alpha, p}}$ is subgaussian with
  \begin{equation}
    \norm[\Big]{
      \norm{ \mu^\epsilon_N - \mu^\epsilon}_{\mathcal W^{-\alpha, p}}
    }_{\psi_2}
    \leq \frac{C d^{1/p} \sqrt{p}}{\sqrt{N}} \norm{\varPhi_\epsilon}_{\mathcal W^{-\alpha, p}}
  \end{equation}
  This immediately implies~\eqref{e:MuNMuSubG}, concluding the proof.
 \end{proof}

It remains to prove Lemmas~\ref{l:subGbound} and~\ref{l:sigmaCalW}, which we do in subsequent sections.

\subsection{The subgaussian norm of \texorpdfstring{$\varPhi_t(x - X_1)$}{} (Lemma \ref{l:subGbound})}\label{s:subGProof}

The key step in the proof of Proposition~\ref{p:empMeasCalW} was the bound on the subgaussian norm of $\varPhi_t(x - X_1)$.
We now present its proof.
\begin{proof}[Proof of Lemma~\ref{l:subGbound}]
  We know (see for instance~\cite[Proposition 2.6.6 (i)]{Vershynin18}) that if for some~$K_1(x,t) > 0$ and all~$\lambda > 0$ we have
  \begin{equation}\label{e:K1bound}
    \P\bigl( \abs{\varPhi_t(x-X_1)} > \lambda \bigr) \leqslant 2 \exp \Bigl(- \frac{\lambda^2}{K_1(x,t)^2}\Bigr)
    \,,
  \end{equation}
  then there exists an absolute constant~$C$ such that
  \begin{equation}\label{e:sigmaK1}
    \sigma(x, t) \defeq \norm{\varPhi_t(x - X_1)}_{\psi_2} \leq C K_1(x, t)
    \,.
  \end{equation}
  Next we note that 
  \begin{equation}
    \lambda \leq K_1(x, t) \sqrt{\ln 2}
    \implies
    2 \exp \Bigl(- \frac{\lambda^2}{K_1(x,t)^2}\Bigr) \geq 1 \geq \P\bigl( \abs{\varPhi_t(x-X_1)} > \lambda \bigr)
    \,. 
  \end{equation}
  Also,
  \begin{equation}
    \lambda \geq (4 \pi t)^{-d/2}
    \implies
    \P\bigl( \abs{\varPhi_t(x-X_1)} > \lambda \bigr) = 0 \leq 2 \exp \Bigl(- \frac{\lambda^2}{K_1(x,t)^2}\Bigr)
    \,,
  \end{equation}
  which again verifies~\eqref{e:K1bound}.
  Thus we only need to ensure that~\eqref{e:K1bound} holds for all $\lambda$ such that
  \begin{equation}\label{e:lambdaInterval}
    K_1 \sqrt{ \ln 2} <  \lambda <  \frac{1}{(4\pi t)^{d/2}}
    \,.
  \end{equation}

  Now define~$r = r(t, \lambda)$ by
  \begin{equation}\label{e:Rdef}
    r \defeq \sqrt{4t \ln \paren[\bigg]{ \frac1{(4\pi t)^{\frac{d}2} \lambda } } }
    \,,
  \end{equation}
  and note
  \begin{equation}
    \P\bigl( \varPhi_t(x-X_1) > \lambda \bigr)
      = \P \paren[\big]{ X_1 \in B(x, r )}
      = \mu \paren[\big]{ B(x, r )}
      \leq V_d M\mu(x) r^d
      \,.
  \end{equation}
  Thus~\eqref{e:K1bound} will hold for all~$\lambda$ satisfying~\eqref{e:lambdaInterval} provided we can ensure
  \begin{equation}\label{e:VdMmu}
    V_d M\mu(x) r^d \leq 2 \exp \Bigl(- \frac{\lambda^2}{K_1(x,t)^2}\Bigr)
    \,.
  \end{equation}
  Clearly this is equivalent to showing that for all
  \begin{equation}\label{e:rinterval}
    0 < r < \sqrt{4t \ln \paren[\bigg]{ \frac1{(4\pi t)^{\frac{d}2} K_1(x, t) (\ln 2)^{1/2} } } }
    \,,
  \end{equation}
  we have
  \begin{equation}\label{e:Dk12ln}
    d K_1^2 { \ln \paren[\Big]{ \frac{\eta(x)}{r} } }\geq \frac{1}{(4 \pi t)^d} e^{-r^2/ 2t}
    \quad\text{where}\quad
    \eta(x) = \paren[\Big]{ \frac{2}{V_d M\mu(x)}}^{1/d}
    \,.
  \end{equation}

  Choose
  \begin{equation}\label{e:K1Def1}
    K_1(x, t) = \frac{1}{(4 \pi t)^{d/2}} e^{-\gamma^2 \eta(x)^2 / 4t}\,,
  \end{equation}
  for some small constant~$\gamma > 0$ that will be chosen shortly.
  Note that both the left and right hand side of~\eqref{e:Dk12ln} are decreasing functions of~$r$.
  Moreover, the inequality~\eqref{e:Dk12ln} clearly holds for~$r = 0$.
  Thus~\eqref{e:Dk12ln} will hold for all~$r$ in the interval~\eqref{e:rinterval}, provided it holds for~$r = r_1$, where
  \begin{equation}
    r_1
      \defeq \sqrt{4t \ln \paren[\bigg]{ \frac1{(4\pi t)^{\frac{d}2} K_1(x, t) (\ln 2)^{1/2} } } }
      = \sqrt{ \gamma^2 \eta(x)^2 + 2 t \ln \paren[\Big]{ \frac{1}{\ln 2}} }
    \,.
  \end{equation}
  When~$r = r_1$, and~$K_1$ is defined by~\eqref{e:K1Def1}, the inequality~\eqref{e:Dk12ln} is equivalent to
  \begin{equation}
    d e^{ - \gamma^2 \eta(x)^2 / 2t} \ln\paren[\Big]{\frac{\eta(x)}{r_1}}
      \geq e^{-r_1^2 / 2t}
      \,.
  \end{equation}
  Through a direct calculation, this in turn is equivalent to having~$\gamma \leq \gamma_1$, where
  \begin{equation}\label{e:gamma1Def}
    \gamma_1 \defeq
	\sqrt{
	  \frac{1}{2^{1/d}} - \frac{2 t \ln \paren[\big]{ \frac{1}{\ln 2}} }{\eta(x)^2}
	}
	\,.
  \end{equation}

  We now divide the analysis into two cases.

  \restartcases
  \case[$\eta(x)^2 > 2^{1 + 1/d} t \abs{\ln \ln 2}$]
  In this case the quantity under the square-root in~\eqref{e:gamma1Def} is positive and we choose~$\gamma = \gamma_1 > 0$.
  Using the definition of~$\eta$ (from~\eqref{e:Dk12ln}) this implies
  \begin{equation}
    K_1(x, t) = \frac{1}{(4 \pi t)^{d/2} \sqrt{\ln 2}}
      \exp\paren[\bigg]{
	\frac{- 2^{3/d}}{4t \paren[\big]{V_d M\mu(x)}^{2/d} }
      }
      \,.
  \end{equation}
  Now the argument above shows~\eqref{e:K1bound} holds for all~$\lambda > 0$, and so using~\eqref{e:sigmaK1} we obtain~\eqref{e:subGbound} as desired.

  \case[$\eta(x)^2 \leq 2^{1 + 1/d} t \abs{\ln \ln 2}$]
  In this case the definition of~$\eta$ (in~\eqref{e:Dk12ln}) implies
  \begin{equation}
    \exp\paren[\bigg]{
      \frac{- 1}{4 t \paren[\big]{ V_d M\mu(x) }^{2/d} }
    }
    \geq \frac{1}{C_1'(d)}
    \quad\text{where}\quad
    C_1'(d) \defeq \exp\paren[\Big]{ \frac{ \abs{\ln \ln 2}}{2^{1 + 1/d}}}
    \,.
  \end{equation}
  Now
  \begin{align}
    \norm{\varPhi_t(x - X_1)}_{\psi_2}
      &\leq C \norm{\varPhi_t(x - X_1)}_{L^\infty}
      \leq \frac{C}{(4 \pi t)^{d/2}}
    \\
    \label{e:subGcaseII}
      &\leq \frac{C C'_1(d)}{(4 \pi t)^{d/2}}
	\exp\paren[\bigg]{
	  \frac{- 1}{ 4 t \paren[\big]{ V_d M\mu(x) }^{2/d} }
	}
      \,,
  \end{align}
  which immediately implies~\eqref{e:subGbound} as desired.
\end{proof}

\subsection{Decay of \texorpdfstring{$\sigma$}{sigma} (Lemma \ref{l:sigmaCalW})}
Given Lemma~\ref{l:subGbound}, one can prove Lemma~\ref{l:sigmaCalW} quickly using the weak~$L^1$ bound on the Hardy--Littlewood maximal function.
We present the proof below.

\begin{proof}[Proof of Lemma~\ref{l:sigmaCalW}]
  Define the function~$F$ by
  \begin{equation}\label{e:Fdef}
    F(\lambda, t)
      = \frac{1}{(4\pi t)^{d/2}} \exp \paren[\bigg]{ \frac{-1}{ 4t \paren{V_d \lambda}^{2/d}} }\,,
      \,.
  \end{equation}
  Using~\eqref{e:subGbound} and the layer-cake formula we see
  \begin{equation}
    \norm{\sigma(\cdot, t)}_{L^p}^p
      \leq C^p \int_{\lambda=0}^\infty \partial_\lambda (F^p) \, \abs[\big]{ \set{M\mu > \lambda} } \, d\lambda
      \leq C^p d \int_{\lambda=0}^\infty \partial_\lambda (F^p) \, \frac{d\lambda}{\lambda}
      \,.
  \end{equation}
  The last inequality above followed from the Stein--Str\"omberg weak~$L^1$ bound for the maximal function (see Theorem~3 in~\cite{SteinStromberg83}).
  Substituting~$r = 1 / (V_d \lambda)^{1/d}$ and using
  \begin{equation}\label{e:drdLambda}
    \frac{d\lambda}{\lambda} = (-d) \frac{dr}{r}
    \qquad
    \partial_\lambda = \frac{-V_d}{d} r^{d+1} \partial_r
  \end{equation}
  we obtain
  \begin{equation}\label{e:sigmatLp}
    \norm{\sigma(\cdot, t)}_{L^p}^p
      \leq -C^p d V_d \int_0^\infty \partial_r (F^p) \, r^d \, dr
      = C^p d^2 V_d \int_0^\infty F^p \, r^{d -1} \, dr
      \,.
  \end{equation}
  The last equality above was obtained by using the fact that~$F(r) r^d \to 0$ as both~$r \to 0$ and~$r \to \infty$, and integrating by parts.

  Next we notice that using radial coordinates to compute~$\norm{\varPhi_t}_{L^p}^p$ gives
  \begin{equation}\label{e:phitLp}
    \norm{\varPhi_t}_{L^p}^p
      = \frac{1}{(4 \pi t)^{pd/2}} \int_{\R^d} \exp\paren[\Big]{ \frac{-p r^2}{4t} }
    A_{d-1} r^{d-1} \, dr
    \,,
  \end{equation}
  where~$A_{d-1}$ is the surface area of the unit sphere in~$\R^d$.
  Using the identity~$d V_d = A_{d-1}$ and~\eqref{e:Fdef} this implies
  \begin{equation}
    \norm{\varPhi_t}_{L^p}^p
      = d V_d \int_0^\infty F^p \, r^{d -1} \, dr
      \,.
  \end{equation}
  Substituting this in~\eqref{e:sigmatLp} implies
  \begin{equation}
    \norm{\sigma(\cdot, t)}_{L^p}^p
      \leq C^p d \norm{\varPhi_t}_{L^p}^p
      \,,
  \end{equation}
  which immediately yields~\eqref{e:intSigmaBd}, concluding the proof.
\end{proof}


\section{Convergence in \texorpdfstring{$\mathscr W^{-\alpha, p}$}{scr W -alpha p}.}\label{s:scrW}

We now prove Theorem~\ref{t:empMeas} when~$\alpha \in \N$.
As before, we do this by proving a slightly more general result.
Let~$\mathscr W^{-\alpha , p} = \mathscr W^{-\alpha, p}(\R^d)$ be the space of all tempered distributions~$u$ for which the norm in~\eqref{e:WMinusAlphaInt} is finite.
We now show that Theorem~\ref{t:empMeas} holds with~$\norm{\cdot}_{W^{-\alpha, p}}$ replaced with~$\norm{\cdot}_{\mathscr W^{-\alpha, p}}$, whether or not~$\alpha \in \N$.
This will immediately imply Theorem~\ref{t:empMeas} for~$\alpha \in \N$, and hence conclude the proof of Theorem~\ref{t:empMeas}.

\begin{proposition}\label{p:empMeasScrW}
  There exists an absolute constant~$C$, independent of~$d$ and~$\mu$, such that for any~$p\in (1,\infty)$ and~$\epsilon \geq 0$ we have the moment bound
  \begin{equation}\label{e:MuNMuScrW}
    \paren[\big]{\E \norm{\mu^\epsilon_N - \mu^\epsilon}_{\mathscr W^{-\alpha, p}}^p}^{1/p}
    \leq
    \frac{C d^{1/p} \sqrt{p}}{\sqrt{N}} \norm{\varPhi_\epsilon}_{\mathscr W^{-\alpha, p}}
    \,,
  \end{equation}
  and the Gaussian tail bound
  \begin{equation}\label{e:MuNMuSubGScr}
    \P\paren[\big]{
    \norm{\mu^\epsilon_N - \mu^\epsilon}_{\mathscr W^{-\alpha, p}} > \lambda
    }
    \leq 2 \exp\paren[\bigg]{ \frac{-N \lambda^2}{C d^{2/p} p \norm{\varPhi_\epsilon}_{\mathscr W^{-\alpha, p}}^2 } }
    \,,
    \quad\text{for every } \lambda > 0\,.
  \end{equation}
  If~$\alpha p > (p - 1) d$ and~$\mu \in \mathscr W^{-\alpha, p}$, then both~\eqref{e:MuNMuScrW} and~\eqref{e:MuNMuSubGScr} also hold for~$\epsilon = 0$.
\end{proposition}

While the proof of Proposition~\ref{p:empMeasScrW} uses the same strategy as that of Proposition~\ref{p:empMeasCalW}, there is one added difficulty.
When working in~$\mathcal W^{-\alpha, p}$, we note~$\E \norm{\mu_N - \mu}_{\mathcal W^{-\alpha, p}}^p$ can be written in terms of~$\E S_N$ as in~\eqref{e:MuNFubini}, after which we use Hoeffding's inequality to proceed.
In the space~$\mathscr W^{-\alpha, p}$ this is no longer possible.
We work around this using a Minkowski type inequality and bound the~$p$-th moment~$\E \norm{\mu_N - \mu}_{\mathscr W^{-\alpha, p}}^p$ in terms of an integral of~$\sigma$.
For clarity of presentation, we first state a lemma bounding~$\norm{\sigma}_{\mathscr W^{-\alpha, p}}$.
\begin{lemma}\label{l:sigmaScrWbd}
  There exists an absolute constant~$C$, independent of~$d$, $\mu$, such that for every~$p \in (1, \infty)$, $\alpha > 0$, $\epsilon \geq 0$ we have
  \begin{equation}\label{e:noSiEpScrW}
    \paren[\bigg]{\int_{\R^d}
      \paren[\Big]{
	\int_0^1 t^{\alpha - 1} \sigma(x, t + \epsilon)^2 \, dt
      }^{p/2} \, dx}^{1/p}
      \leq C d^{1/p} \norm{\varPhi_\epsilon}_{\mathscr W^{-\alpha, p}}
      \,.
  \end{equation}
\end{lemma}
\begin{remark}
  As with Remark~\ref{r:alphaBigCalW}, note that when~$\alpha > d/q$ (i.e.\ when~$\alpha p > (p-1) d$), the right hand side of~\eqref{e:noSiEpScrW} is finite when~$\epsilon = 0$.
\end{remark}
Momentarily postponing the proof of Lemma~\ref{l:sigmaScrWbd}, we now prove Proposition~\ref{p:empMeasScrW}.

\begin{proof}[Proof of Proposition~\ref{p:empMeasScrW}]
  Notice
  \begin{equation}\label{e:MuNepMinusMuSrcW}
    \E \norm{\mu^\epsilon_N - \mu}_{\mathscr W^{-\alpha, p}}^p
      = \int_{\R^d} \E J_{N, \alpha, \epsilon}(x)^p \, dx
  \end{equation}
  where~$J_{N, \alpha, \epsilon}(x)$ is the random variable
  \begin{equation}
    J_{N, \alpha, \epsilon}(x) \defeq  \paren[\bigg]{
      \int_0^1 t^{\alpha - 1} S_N(x, t+\epsilon)^2 \, dt
    }^{1/2}
  \end{equation}
  and~$S_N$ is defined in~\eqref{e:SNdef}.
  Now let~$\norm{\cdot}_{\psi_1}$ denote the subexponential norm, and recall (see for instance~\cite{Vershynin18}, Lemma 2.7.6) that for any random variable~$Z$ we have
  \begin{equation}
    \norm{Z^2}_{\psi_1} = \norm{Z}_{\psi_2}^2
    \,.
  \end{equation}
  Since~$\norm{\cdot}_{\psi_1}$ is a norm, the above identity and Minkowski's inequality imply
  \begin{align}
    \norm{J_{N, \alpha, \epsilon}(x)}_{\psi_2}
      &= \norm{J_{N, \alpha, \epsilon}(x)^2}_{\psi_1}^{1/2}
      \leq
	\paren[\bigg]{
	  \int_0^1 t^{\alpha - 1} \norm{S_N(x, t+\epsilon)^2}_{\psi_1} \, dt
	}^{1/2}
    \\
      &= \paren[\bigg]{
	  \int_0^1 t^{\alpha - 1} \norm{S_N(x, t+\epsilon)}_{\psi_2}^2 \, dt
	}^{1/2}
    \\
      \label{e:InPsi2}
      &\overset{\mathclap{\eqref{e:SnPsi2}}}\leq \frac{C}{\sqrt{N}}
	\paren[\bigg]{
	  \int_0^1 t^{\alpha - 1} \sigma(x, t+\epsilon)^2 \, dt
	}^{1/2}
	\,.
  \end{align}

  Substituting the above in~\eqref{e:MuNepMinusMuSrcW} implies
  \begin{align}
    \paren[\big]{\E \norm{\mu^\epsilon_N - \mu}_{\mathscr W^{-\alpha, p}}^p}^{1/p}
      &\leq C \sqrt{p} \norm[\Big]{ \norm{J_{N, \alpha, \epsilon}}_{\psi_2} }_{L^p}
      %
      \\
      &\overset{\mathclap{\eqref{e:InPsi2}}}{\leq}~
	\frac{C \sqrt{p}}{\sqrt{N}}
	\paren[\bigg]{
	  \int_{\R^d}
	    \paren[\bigg]{
	      \int_0^1 t^{\alpha - 1} \sigma(x, t+\epsilon)^2 \, dt
	    }^{p/2} \, dx
	}^{1/p}
	\,.
  \end{align}
  Combined with Lemma~\ref{l:sigmaScrWbd} this yields~\eqref{e:MuNMuScrW} as desired.
  \medskip

Next we turn to the Gaussian tail bound~\eqref{e:MuNMuSubGScr}.
  The argument is almost identical to the proof of~\eqref{e:MuNMuSubG}.
  As with the proof of~\eqref{e:MuNMuSubG}, the starting point is to notice that 
\begin{equation}
  \|\mu^\e_N-\mu^\e\|_{\mathscr{W}^{-\alpha,p}} = f_N(X_1,\cdots, X_N),
\end{equation}
where 
\begin{equation}
  f_N(y_1,\cdots,y_N) \defeq \|g_N (\cdot, \cdot; y_1,\cdots,y_N)\|_{\star}\,.
\end{equation}
Here~$g_N$ is the same function defined in~\eqref{e:gNdef}, and~$\norm{\cdot}_\star$ is the weighted norm
\begin{equation}
  \|h\|_{\star} \defeq \biggl( \int_{\Rd} \biggl(\int_0^1 t^{\alpha-1} h(x,t)^2\,dt \biggr)^{\frac{p}2}\,dx \biggr)^{\frac1{p}} \,. 
\end{equation}

We now claim for any~$i \in \set{1, \dots, N}$, $y_i, y_i' \in \R^d$ the function~$f_N$ satisfies
\begin{equation}\label{e:bddDiff1}
  \abs{f_N(y_1, \dots, y_N) - f_N(y_1, \dots, y_{i-1}, y_i', y_{i+1}, \dots, y_N)}
    \leq \frac{2}{N} \norm{\varPhi_\epsilon}_{\mathscr W^{-\alpha, p}}
    \,.
\end{equation}
The proof of~\eqref{e:bddDiff} is very similar to that of~\eqref{e:bddDiff}.
Indeed, by symmetry it suffices to prove~\eqref{e:bddDiff1} for~$i = 1$.
For this choose any~$y_1,y_1'\in\Rd$ and~$\tilde{y}= (y_2,\cdots, y_N) \in (\Rd)^{N-1}$.
The reverse triangle inequality implies
\begin{equation}
	|f_N(y_1,\tilde{y}) - f_N(y_1',\tilde{y})| \leqslant \|h(\cdot,\cdot; y_1,y_1') \|_{\star}\,,
\end{equation}
where~$h$ is the same function defined in~\eqref{e:hNDef}.
This immediately implies
  \begin{equation}
    \norm{h(\cdot, \cdot; y_1, y_1')}_{\star}
      \leq \frac{2}{N} \norm{\varPhi_\epsilon}_{\mathscr W^{-\alpha, p}}
      \,,
  \end{equation}
from which~\eqref{e:bddDiff1} follows.

Given~\eqref{e:bddDiff1} we obtain~\eqref{e:MuNMuSubGScr} using McDiarmid's inequality and the moment bound~\eqref{e:MuNMuScrW} exactly as in the proof of~\eqref{e:MuNMuSubG}.
\end{proof}

It remains to prove Lemma~\ref{l:sigmaScrWbd}.
The proof is a direct consequence of Lemma~\ref{l:subGbound} and the Stein--Str\"omberg weak~$L^1$ bound on the Hardy--Littlewood maximal function (Theorem~3 in~\cite{SteinStromberg83}).
\begin{proof}[Proof of Lemma~\ref{l:sigmaScrWbd}]
Define the function~$F$ by
\begin{equation*}
  F(\lambda) \defeq \biggl(\int_0^1 t^{\alpha-1}  \frac{1}{(4\pi (t+\e)^d)} \exp  \biggl(\frac{-1}{ 2(t+\e) (V_d\lambda)^{2/d}}\biggr)\,dt\biggr)^{\frac{p}2}\,,
\end{equation*}
and observe
\begin{align}
  \norm{\sigma_\epsilon}_{\mathscr W^{-\alpha, p}}^p
  &\overset{\mathclap{\eqref{e:subGbound}}}\leq~ C^p \int_{\Rd} F(M\mu(x)) \,dx
    = C^p \int_0^\infty \partial_\lambda F(\lambda)\bigl|\{ M\mu(x) > \lambda \}\bigr|\,d\lambda
  \\
  \label{e:siEpScrW}
  &\leqslant C^p d \int_0^\infty  \partial_\lambda F\,  \frac{d\lambda}{\lambda}
  \,.
\end{align}
  Making the same change of variables~$r = 1/(V_d \lambda)^{1/d}$, and using~\eqref{e:drdLambda} shows
  \begin{equation}\label{e:intDLaF1}
    \int_0^\infty  \partial_\lambda F(\lambda) \,  \frac{d\lambda}{\lambda}
      = -V_d \int_0^\infty  \partial_r F(r) \,  r^d \, dr
      = d V_d \int_0^\infty  F(r)\,  r^{d-1} \, dr
      \,.
  \end{equation}
  To obtain the last equality above we integrated by parts and used the fact that~$r^d F(r)$ vanishes at both~$0$ and~$\infty$.
  Computing~$\norm{\varPhi_\epsilon}$ using radial coordinates and using the identity~$d V_d =A_{d-1}$ we see
  \begin{equation}
    d V_d \int_0^\infty  F(r)\,  r^{d-1} \, dr
      = \norm{\varPhi_\epsilon}_{\mathscr W^{-\alpha, p}}^p
  \end{equation}
  Combining this with~\eqref{e:siEpScrW} and~\eqref{e:intDLaF1} yields~\eqref{e:noSiEpScrW} as desired.
\end{proof}

\section{Exact second moments in \texorpdfstring{$H^{-\alpha}$}{H -alpha} (Proposition \ref{p:EmpMeasureHAlpha}).}\label{s:Halpha}

The proof of Proposition~\ref{p:EmpMeasureHAlpha} is a short and elementary calculation.
\begin{proof}[Proof of Proposition~\ref{p:EmpMeasureHAlpha}]
  For any~$\epsilon > 0$ and~$N \in \N$ we have
  \begin{equation*}
    \| \mu^\epsilon_N - \mu^\epsilon \|_{H^{-\alpha}}^2
      = \int_0^1 \int_{\Rd} t^{\alpha - 1} \paren[\Big]{
	  \frac{1}{N} \sum_{i = 1}^N \varPhi_{t + \epsilon}(x - X_i) - \mu_{t+\epsilon}(x)
	}^2
	\,dx \, dt
      \,.
  \end{equation*}
  Notice that for every~$x$ and $t$ the summands~$\varPhi_{t+\epsilon}(x - X_i) - \mu_{t+\epsilon}(x)$ are centered i.i.d.\ random variables.
  Thus
  \begin{equation}\label{e:MuNMuHAlpha}
    \E\| \mu^\epsilon_N - \mu^\epsilon \|_{H^{-\alpha}}^2
      = \frac{1}{N} \int_0^1 \int_{\Rd} t^{\alpha - 1} \paren[\big]{
	  \E \varPhi_{t + \epsilon}(x - X_1)^2 - \mu_{t+\epsilon}(x)^2
	}
	\,dx \, dt
      \,.
  \end{equation}

  Notice
  \begin{align}
    \MoveEqLeft
      \int_0^1 \int_{\Rd} t^{\alpha - 1} 
	  \E \varPhi_{t + \epsilon}(x - X_1)^2
	\,dx \, dt
      = \int_0^1 \int_{\Rd \times \R^d} t^{\alpha - 1} 
	  \varPhi_{t + \epsilon}(x - y)^2
	  \, d\mu(y)
	\,dx \, dt
    \\
      &= \int_0^1 \int_{\Rd \times \R^d} t^{\alpha - 1} 
	  \varPhi_{t + \epsilon}(z)^2
	  \, dz \, d\mu(y) \, dt
      = \int_{\R^d} \norm{\varPhi_\epsilon}_{H^{-\alpha}}^2 \, d\mu(y)
      = \norm{\varPhi_\epsilon}_{H^{-\alpha}}^2
      \,.
  \end{align}
  Also, by definition,
  \begin{equation}
    \int_0^1 \int_{\Rd} t^{\alpha - 1}
	\mu_{t+\epsilon}(x)^2
      \,dx \, dt
    = \norm{\mu^\epsilon}_{H^{-\alpha}}^2
    \,.
  \end{equation}
  Thus~\eqref{e:MuNMuHAlpha} reduces to~\eqref{e:MuNMinusMuEpHAlpha} as desired.
  When~$\alpha > d/2$, all terms in~\eqref{e:MuNMuHAlpha} remain finite when~$\epsilon = 0$, and in this case~\eqref{e:MuNMuHAlpha} also holds when~$\epsilon = 0$.

  The proof of~\eqref{e:MuNMuSubGHalpha} is almost identical to that of~\eqref{e:MuNMuSubG} and~\eqref{e:MuNMuSubGScr}.
  In the proofs of both~\eqref{e:MuNMuSubG} and~\eqref{e:MuNMuSubGScr} we note that the dimensional factor~$d^{1/p}$ only arises through the bound on the mean (inequalities~\eqref{e:MuNMuCalW} and~\eqref{e:MuNMuScrW} respectively).
  In our case the mean is the exactly given by~\eqref{e:MuNMinusMuEpHAlpha}, which doesn't have the~$d^{1/p}$ factor.
  As a result this factor is also absent from the Gaussian tail bound~\eqref{e:MuNMuSubGHalpha}.
\end{proof}

\appendix
\section{The \texorpdfstring{$\mathcal W^{-\alpha, p}$}{cal W} and \texorpdfstring{$\mathscr W^{-\alpha, p}$}{scr W} norms of Gaussians.}

Since the bounds in Propositions~\ref{p:empMeasCalW} and~\ref{p:empMeasScrW} involve~$\norm{\varPhi_\epsilon}_{\mathcal W^{-\alpha, p}}$ and~$\norm{\varPhi_\epsilon}_{\mathscr W^{-\alpha, p}}$, we now compute them explicitly.
A scaling argument quickly shows their behavior in~$\epsilon$ as~$\epsilon \to 0$.
However, their explicit dimensional dependence is useful and isn't easily available in the existing literature.
We present it here for reference.

\begin{proposition}\label{p:PhiEpBounds}
  For every~$\epsilon > 0$ we have
  \begin{align}
    \label{e:PhiEpCalW}
    \|\varPhi_\e\|_{\mathcal{W}^{-\alpha,p}}
      &= \frac{\e^{\frac{1}{2}(\alpha - \frac{d}{q})}}{p^{\frac{d}{2p}} (4\pi)^{\frac{d}{2q}}}  \mathcal B_0(\e,\alpha,p,d)^{1/p}\,,
  \\
    \label{e:PhiEpScrW}
    \|\varPhi_\e\|_{\mathscr{W}^{-\alpha,p}}
      &= \frac{\e^{\frac{1}{2}(\alpha - \frac{d}{q})}}{(4\pi)^{\frac{d}2}} \mathscr B_0(\epsilon, \alpha, p, d)^{1/p}\,,
  \end{align}
  where~$\mathcal B_0$ and~$\mathscr B_0$ are given by
    \begin{align}\label{e:B0Def}
      \mathcal B_0(\epsilon)
      &= \mathcal B_0(\epsilon, \alpha, p, d)
      \defeq
      \int_0^{1/\epsilon} \frac{s^{\frac{\alpha p}{2} - 1}}{(1+s)^{(p-1)\frac{d}{2}}} \, ds
      \,,
    \\
      \label{e:C0Def}
      \mathscr B_0(\epsilon) &= \mathscr B_0(\epsilon,\alpha,p,d) = \int_{\Rd} \biggl(\int_0^{\frac1{\e}} \frac{s^{\alpha-1}}{ (1+s)^d} \exp \biggl(- \frac{|y|^2}{2(1+s)}\biggr)\,ds\biggr)^{\frac{p}2}\,dy
      \,.
    \end{align}
  \begin{enumerate}
    \item
      When~$\alpha > d/q$
      \begin{align}
	\label{e:calBalphaBig}
	\lim_{\epsilon \to 0}\e^{\frac{1}{2}(\alpha - \frac{d}{q})} \mathcal B_0(\epsilon)^{1/p}
	  &= \Bigl(\frac{2}{\alpha p - d(p-1)} \Bigr)^{\frac1{p}}\,,\\
	\label{e:calBalphaSmall}
	\lim_{\epsilon \to 0}\e^{\frac{1}{2}(\alpha - \frac{d}{q})} \mathscr B_0(\epsilon)^{1/p}
	  &\in (0, \infty)
	\,,
      \end{align}
      and~$\norm{\delta_0}_{\mathcal W^{-\alpha, p}}$, $\norm{\delta_0}_{\mathscr W^{-\alpha, p}}$ are the limits as~$\epsilon \to 0$ of the expressions in~\eqref{e:PhiEpCalW} and~\eqref{e:PhiEpScrW} respectively.

    \item
      When~$\alpha = d/q$
      \begin{equation}
	\lim_{\epsilon \to 0} \paren[\bigg]{\frac{\mathcal B_0(\epsilon)}{\abs{\ln \epsilon}}}^{1/p}
	  = 1\,,
	\quad\text{and}\quad
0<	\lim_{\epsilon \to 0} \paren[\bigg]{\frac{\mathscr B_0(\epsilon)}{\abs{\ln \epsilon}}}^{1/p}
	  \leqslant A_{d-1}^{\frac1{p}} \Gamma \Bigl( \frac{d}{p}\Bigr)^{\frac1{p}}
	\,.
      \end{equation}

    \item
      When~$0 < \alpha < d/q$, both~$\mathcal B_0(\epsilon)$ and~$\mathscr B_0(\epsilon)$ converge to finite, non-zero constants as~$\epsilon \to 0$.
      Hence both norms~$\norm{\varPhi_\epsilon}_{\mathcal W^{-\alpha, p}}$ and~$\norm{\varPhi_\epsilon}_{\mathscr W^{-\alpha, p}}$ diverge like~$1/\epsilon^{\frac{1}{2}(\frac{d}{q} - \alpha)}$ as~$\epsilon \to 0$.
  \end{enumerate}
\end{proposition}

For clarity of presentation we divide the proof into two parts, each addressing bounds in~$\mathcal W^{-\alpha, p}$ and~$\mathscr W^{-\alpha, p}$ respectively.
\subsection{Computation of \texorpdfstring{$\norm{\varPhi_\epsilon}_{\mathcal W^{-\alpha, p}}$}{norm(varPhi) cal W -alpha p}}

\begin{proof}[Proof of Proposition~\ref{p:PhiEpBounds} for~$\mathcal W^{-\alpha, p}$]
Using the semigroup property of the heat kernel we compute
\begin{align*}
&\|\varPhi_\e\|_{\mathcal{W}^{-\alpha,p}}^p = \int_0^1 \int_{\Rd} t^{\frac{\alpha p}2 - 1}   \varPhi_{t+\e}(x)^p\,dx\,dt \\
&\quad = \int_0^1 \int_{\Rd} \frac{t^{\frac{\alpha p}2 - 1} }{(4\pi (t+\e))^{\frac{dp}2}} \exp \biggl( - \frac{p|x|^2}{4(t+\e)}\biggr)\,dx\,dt \\
&\quad = \int_0^1  \frac{t^{\frac{\alpha p}2 - 1} }{(4\pi (t+\e))^{\frac{dp}2}} \Bigl(\frac{4\pi (t+\e)}{p}\Bigr)^{\frac{d}2} \\
&\quad = \frac1{p^{d/2} (4\pi)^{\frac{d}2(p-1)}} \int_0^1 \frac{t^{\frac{\alpha p}2 - 1}}{(t+\e)^{\frac{d(p-1)}2}}\,dt = \frac{\e^{\frac{\alpha p}2 - \frac{d(p-1)}{2} }}{p^{d/2} (4\pi)^{\frac{d}2(p-1)}} \mathcal{B}_0(\e,\alpha,p,d)\,,
\end{align*}
where the last equality follows from a change of variables.
Taking the~$p-$th root yields~\eqref{e:PhiEpCalW}. 

We now prove items (1)--(3) for the quantity~$\mathcal{B}_0(\e)$.
First, when~$\alpha > d/q$, the dominated convergence theorem yields
\begin{equation}
	  \lim_{\e\to 0^+} \biggl( \e^{\frac{\alpha p}2 - \frac{d(p-1)}{2} } \mathcal{B}_0(\e,\alpha,p,d)  \biggr)^{1/p}
	= \Bigl(\frac{2}{\alpha p - d(p-1)} \Bigr)^{\frac1{p}}\,. 
\end{equation}
  This proves~\eqref{e:calBalphaBig} as desired.
  Moreover, in this regime, the Sobolev embedding theorem implies~$\varPhi_\epsilon \to \delta_0$ in  $W^{-\alpha, p}$ as~$\epsilon \to 0$ and hence
\begin{equation*}
	\|\delta_0\|_{\mathcal{W}^{-\alpha,p}} = 	\lim_{\e\to 0^+} \|\varPhi_\e\|_{\mathcal{W}^{-\alpha,p}}
	\,,
\end{equation*}
as claimed.
\smallskip 

Turning to item~(2) for the quantity~$\mathcal{B}_0(\e)$, we compute directly by the l'H\^opital rule that 
\begin{align}
  \lim_{\e \to 0^+} \frac{\mathcal{B}_0(\e)}{|\ln \e|} &= - \lim_{\e \to 0^+} \frac1{\ln \e} \int_0^{1/\e} \frac{s^{\frac{\alpha p}2 - 1}}{(1 + s)^{(p-1)\frac{d}2}}  \,ds \\
	&= - \lim_{\e \to 0^+} \e \frac{(1/\e)^{\frac{\alpha p }2 - 1}}{(1 + 1/\e)^{(p-1)\frac{d}2}} (- \frac1{\e^2}) = \lim_{\e \to 0^+} \e^{\frac{(p-1)d - \alpha p}2} = 1\,. 
\end{align}
\smallskip 
Finally, to prove item~(3), we now assume~$0 < \alpha < d/q$.
  In this case
\begin{align}
  \mathcal{B}_0 (\e)  &\leqslant \int_0^1 s^{\frac{\alpha p}2 - 1} \,ds + \int_1^{\frac1{\e}} s^{\frac{\alpha p}2 - (p-1)\frac{d}2 - 1} \,ds \\
	&\leqslant \frac2{\alpha p} + \frac2{|\alpha p - (p-1)d|} \bigl( 1 - \e^{(p-1)\frac{d}2 - \frac{\alpha p}2}\bigr) \,, 
\end{align}
and so monotone convergence theorem completes item (3) for the quantity~$\mathcal{B}_0(\e)$, yielding that 
\begin{equation*}
	\lim_{\e \to 0^+} \mathcal{B}_0(\e) = \int_0^\infty \frac{s^{\frac{\alpha p}2 - 1} }{(1+s)^{\frac{(p-1)d}2}} \in (0,\infty)\,. 
\end{equation*}
\end{proof}

\subsection{Computation of \texorpdfstring{$\norm{\varPhi_\epsilon}_{\mathscr W^{-\alpha, p}}$}{norm(varPhi) scr W -alpha p}}
Before we can estimate~$\norm{\varPhi}_{\mathscr W^{-\alpha, p}}$, we need a calculus lemma concerning asymptotics of the inner integral in~\eqref{e:WMinusAlphaInt}.
\begin{lemma}
	\label{l.Iasymptotics}
	For any~$r \geq 0$ let 
	\begin{equation*}
		I_\e(r):= \int_0^{\frac1{\e}} \frac{s^{\alpha - 1}}{(1+s)^d} \exp \biggl( - \frac{r^2}{2(1+s)} \biggr)\,ds\,.  
	\end{equation*}
For all~$r \geq 0$, then~$I_\e(r)$ satisfies the bounds
\begin{equation}
\label{e.unify}
I_\e(r) \leqslant \begin{dcases}
	\frac1{\alpha} + \frac1{\alpha - d} \Bigl( \frac1{\e^{\alpha -d }} - 1\Bigr) &  \alpha \neq d\,\\
		\frac1{\alpha} + \log \frac1{\e}  &  \alpha = d\,.
\end{dcases}
\end{equation}
Additionally, for every~$r > 0$, we have 
\begin{equation} \label{e.decayiny}
  \frac{C_\alpha}{ r^{2(d- \alpha)}} \int_{\frac{\e r^2}{2(1+\e)}}^{\frac{r^2}2} t^{d-\alpha - 1} e^{-t}\,dt
    \leq I_\e(r)
    \leqslant 	 \frac1{r^{2(d-\alpha)}} \int_{ \frac{\e r^2}{2(1+\e)} }^{\frac{r^2}{2}} t^{d-\alpha- 1} e^{-t}\,dt\,. 
\end{equation}
where~$C_\alpha \defeq \min (1 , \frac1{2^{\alpha-1}})$.  
\end{lemma}

\begin{proof}
  We first note the elementary bound
	\begin{equation} \label{e.y<=1alphaneqd}
	I_\e(r) \leqslant \int_0^{\frac1{\e}}  \frac{s^{\alpha - 1}}{(1+s)^d}\,ds \leqslant \int_0^1  s^{\alpha - 1}\,ds + \int_1^{\frac1{\e}} s^{\alpha - d - 1}\,ds = \frac1{\alpha} + \frac1{\alpha - d} \Bigl( \frac1{\e^{\alpha -d }} - 1\Bigr)\,,
	\end{equation}
if~$\alpha\neq d$, and
	\begin{equation} \label{e.y<=1alphaeqd}
	I_\e(r) \leqslant \int_0^{\frac1{\e}}  \frac{s^{\alpha - 1}}{(1+s)^d}\,ds \leqslant \int_0^1  s^{\alpha - 1}\,ds + \int_1^{\frac1{\e}} s^{ - 1}\,ds = \frac1{\alpha} + \log \frac1{\e}\,. 
\end{equation}
when~$\alpha=d$. This is the proof of~\eqref{e.unify}. 

We also need a different bound that captures decay in~$r$ for~$r$ large. For this,  the change of variables~$\frac{r^2}{2(1+s)} = t$, gives
\begin{equation} \label{e.Iepscov}
	I_\e(r) = \int_{\frac{\e r^2}{2(1+\e)} }^{\frac{r^2}{2}}  \biggl( \frac{r^2}{t} - 1\biggr)^{\alpha-1} \frac{t^{d-1}}{r^{2(d-1)}}  e^{-t} \frac{\,dt}{t}\,. 
\end{equation}
Since~$(r^2 - t)^{\alpha -1} \leqslant r^{2(\alpha-1)}$ on the interval of integration, we obtain 
\begin{equation*}
	I_\e(r) \leqslant \frac1{r^{2(d-\alpha)}} \int_{ \frac{\e r^2}{2(1+\e)} }^{\frac{r^2}2} t^{d-\alpha- 1} e^{-t}\,dt  \,.
\end{equation*}
This yields the upper bound in~\eqref{e.decayiny}.  

For the lower bound, using~\eqref{e.Iepscov} one more time and that on the interval of integration, we have~$r^2 \geqslant r^2 \bigl( 1 - \frac{\e}{2(1+\e)}\bigr) \geqslant r^2 - t \geqslant \frac{r^2}{2}$, we compute
\begin{equation*}
	I_\e(r) \geqslant \frac{C_\alpha}{ r^{2(d- \alpha)}} \int_{\frac{\e r^2}{2(1+\e)}}^{\frac{r^2}2} t^{d-\alpha - 1} e^{-t}\,dt\,,
\end{equation*}
with~$C_\alpha = 1$ if~$\alpha \leqslant 1$ and~$C_\alpha = \frac1{2^{\alpha-1}}$ if~$\alpha \geqslant 1$.
  This yields the lower bound in~\eqref{e.decayiny}

\end{proof}
%
%
With these observations about~$I_\e$, we are ready to prove Lemma~\ref{p:PhiEpBounds} for~$\mathscr W^{-\alpha, p}$. 

\begin{proof}[Proof of Proposition~\ref{p:PhiEpBounds} for~$\mathscr W^{-\alpha, p}$]
We compute 
\begin{align}
  \|\varPhi_\e\|_{\mathscr{W}^{-\alpha,p}}^p &= \int_{\Rd} \biggl(\int_0^1 t^{\alpha-1} \varPhi_{t+\e}(x)^2\,dt\biggr)^{\frac{p}2}\,dx\\
	&= \int_{\Rd} \biggl(\int_0^1 \frac{t^{\alpha-1}}{(4\pi(t+\e))^d} \exp\biggl(- \frac{|x|^2}{2(t+\e)}\biggr)\,dt\biggr)^{\frac{p}2}\,dx
\end{align}
Here, make the change of variables~$t = \e s, y = \frac{x}{\sqrt{\e}}$, to obtain 
\begin{align}
  \|\varPhi_\e\|_{\mathscr{W}^{-\alpha,p}}^p  &= \e^{\frac{\alpha p}2 - (p - 1)\frac{d}2 }\int_{\Rd} \biggl(\int_0^{\frac1{\e}} \frac{s^{\alpha-1}}{(4\pi (1+s))^d} \exp \biggl(- \frac{|y|^2}{2(1+s)}\biggr)\,ds\biggr)^{\frac{p}2}\,dy\\
	&\defeq \e^{\frac{\alpha p}2 - (p - 1)\frac{d}2 } \mathscr B_0(\e,\alpha,d,p)\,.
\end{align}
Upon taking the~$p$-th root of both sides, we get~\eqref{e:PhiEpScrW}. We turn to analyzing  the~$\e \to 0^+$ behavior of~$\mathscr B_0$, thereby verifying items (1)-(3) in the different regimes. 

 We must divide our analysis in a few cases.

\restartcases
  \case[$\alpha > d> \paren{p-1}d/p$]
  Lemma~\ref{l.Iasymptotics} show that for small~$\e$ we have
 \begin{align}
\mathscr{B}_0(\e,\alpha,d,p)
   &= 	A_{d-1}\ \int_0^\infty I_\e(r)^{\frac{p}2} r^{d-1}\,dr
   \\
   \label{e:supcrit2}
   &= A_{d-1}  \biggl( \int_0^{\sqrt{\e}} + \int_{\sqrt{\e}}^\infty\biggr) I_\e(r)^{\frac{p}2} r^{d-1}\,dr
   = A_{d-1}\bigl(J_1 + J_2\bigr)
   \,.
  \end{align} 

  Using~\eqref{e.unify} we see
  \begin{equation}
    J_1 \leq    \biggl( \int_0^{\sqrt{\e}} \frac2{(\alpha - d)^{\frac{p}2}} \e^{(d-\alpha)p/2} r^{d-1}\,dr \biggr) = \frac{2}{(\alpha-d)^{p/2}} \e^{\frac12 ((p-1)d - \alpha p)} 
    \,.
  \end{equation}
  For~$J_2$ we use~\eqref{e.decayiny} to obtain
  \begin{equation}
    J_2 = \int_{\sqrt{\e}}^\infty
      I_\e(r)^{\frac{p}2} r^{d-1}\,dr
      \leq \int_{\sqrt{\epsilon}}^\infty r^{\alpha p - (p-1)d- 1}
	\paren[\bigg]{
	  \int_{ \frac{\e r^2}{2(1+\e)} }^{\frac{r^2}{2}} t^{d-\alpha- 1} e^{-t}\,dt
	}^{p/2}\,dr 
	\,.
  \end{equation}
  Making the change of variables~$s = \sqrt{\e} r$ in the last integral, we obtain 
  \begin{align}
  	J_2 &\leqslant  \e^{\frac12 ( (p-1)d - \alpha p )} \int_1^\infty s^{\alpha p - (p-1)d -1}  \biggl( \int_{\frac{s^2}{4}}^{\infty} t^{d-\alpha-1} e^{-t}\,dt \biggr)^{p/2}\,ds\\
  	&\leqslant \e^{\frac12 ( (p-1)d - \alpha p )} \int_1^\infty s^{\alpha p - (p-1)d -1}  e^{-\frac{p s^2}{8}} \frac1{(\alpha - d)^{\frac{p}2}}s^{p(d-\alpha)}\,ds \\
  	& =  \frac{\e^{\frac12 ( (p-1)d - \alpha p )}}{(\alpha - d)^{p/2}} \int_1^\infty s^{d-1} e^{-\frac{p s^2}{8}} \,ds \,. 
  \end{align}

Combining upper bounds for~$J_1$ and~$J_2$ we find
\begin{equation}\label{e.uppersupercritical}
  \mathscr{B}_0(\e,\alpha,d,p)
    \leq \frac{A_{d-1}}{(\alpha - d)^{p/2}}\biggl(2 + \int_1^\infty s^{d-1}e^{-\frac{ps^2}{8}}\,ds \biggr) \e^{\frac{d(p-1) -\alpha p}{2}}  \,,
 \end{equation}
so that in this regimes, the gaussians~$\{\varPhi_\e\}_{\e>0}$ have uniformly bounded~$\mathscr{W}^{-\alpha,p}$ norms. 

In order to complete the proof of~\eqref{e:calBalphaSmall}, we must show a matching lower bound on~$\mathscr{B}_0(\e)$. For this,  using the lower bound portion of~\eqref{e.decayiny} we compute, for~$\e$ small enough (so that~$e^{-\frac{p \e^2}{4}} \geqslant \frac12$),
\begin{align*}
&	\mathscr{B}_0(\e) \geqslant \int_{|y|\leqslant \sqrt{\e}} I_\e(|y|)^{\frac{p}2} \,dy \geqslant C_\alpha A_{d-1}\int_{0}^{\sqrt{\e}} \frac{r^{d-1}}{r^{(d-\alpha)p}} \biggl(\int_{\frac{\e r^2}{2(1+\e)}}^{\frac{r^2}2} t^{d-\alpha-1}e^{-t}\,dt\biggr)^{\frac{p}2}\,dr\\ &\quad 
	\geqslant \frac{A_{d-1}}{(\alpha - d)^p}\int_0^{\sqrt{\e}} \frac{r^{d-1}}{r^{(d-\alpha)p}} e^{- \frac{p r^2}{4}} \biggl( \Bigl( \frac{\e r^2}{2(1+\e)}  \Bigr)^{d-\alpha} - \Bigl( \frac{r^2}2\Bigr)^{d-\alpha} \biggr)^{\frac{p}2} \,dr\\ &\quad 
	\geqslant \frac{A_{d-1}}{2^{1 + \frac{(d-\alpha)p}2}(\alpha - d)^p} \int_0^{\sqrt{\e}} r^{d-1} \frac{\e^{\frac{(d-\alpha)p}2}}2 \,dr =  \frac{A_{d-1}}{2^{2 + \frac{(d-\alpha)p}2}(\alpha - d)^pd}\e^{\frac{d(p-1) -\alpha p}{2}} \,,
\end{align*}
which matches the scaling of the upper bound from~\eqref{e.uppersupercritical}. Combining the previous two displays completes the proof of~\eqref{e:calBalphaSmall}.

  \case[$\frac{p-1}{p}d \leqslant \alpha \leqslant d$]
In this case, we compute 
\begin{align} \label{e.computation-sub}
&\mathscr{B}_0(\e,\alpha,d,p)  \\
&\quad = 	A_{d-1} \int_0^\infty I_\e(r)^{\frac{p}2} r^{d-1}\,dr \\
&\quad = A_{d-1}  \biggl( \int_0^{\sqrt{\e}} + \int_{\sqrt{\e}}^{\frac1{\sqrt{\e}}} + \int_{\frac1{\sqrt{\e}}}^\infty\biggr) I_\e(r)^{\frac{p}2} r^{d-1}\,dr \,. 
\end{align}
Let us analyze each of the pieces in turn, starting first with the middle piece (which will turn out to be most important). To analyze this middle term we must observe that  when~$d>\alpha$, for~$r>0,$ 
 then
\begin{equation} \label{e.gammabound}
	I_\e(r) \leqslant \frac1{r^{2(d-\alpha)}} \int_0^\infty t^{d-\alpha -1} e^{-t}\,dt = \frac{\Gamma(d-\alpha)}{r^{2(d-\alpha)}}\,. 
\end{equation}
Inserting~\eqref{e.gammabound} we obtain
\begin{align}
  \MoveEqLeft
\int_{\sqrt{\e}}^{\frac1{\sqrt{\e}}} I_\e(r)^{\frac{p}2} r^{d-1} \,dr
  \leqslant  \Gamma(d-\alpha)^{\frac{p}2} \int_{\sqrt{\e}}^{\frac1{\sqrt{\e}}}  \frac1{r^{p(d-\alpha)}} r^{d-1}\,dr
  \\
 \label{e.middlepiece}
  &\leqslant
    \Gamma(d-\alpha)^{\frac{p}2}  \times \begin{dcases}
 \ln \frac1{\e} & \mbox{ if } \alpha p = d(p-1)\\
\frac1{|\alpha p - (p-1)d|} \e^{\frac{ d(p-1)- \alpha p}2 }  & \mbox{ if } d(p-1) < \alpha p \leqslant dp\,. 
\end{dcases}
\end{align}
Now we turning to the piece on~$[0,\sqrt{\e}]$ we compute directly that when~$\alpha \leqslant d$, by inserting the uniform bound~\eqref{e.unify} that  
\begin{equation} \label{e.firstpiece}
	\int_{0}^{\sqrt{\e}}I_{\e}(r)^{\frac{p}2} r^{d-1}\,dr \leqslant C(\alpha,d,p) \e^{\frac{d}2}\bigl( 1 + (\log^{\frac{p}2} \frac1{\e}) 1_{\alpha = d} \bigr)\,,
\end{equation}
where the notation~$1_{\alpha = d}$ serves to indicate that this logarithmic correction is only present when~$\alpha=d$. 

Finally, the contribution of the piece~$\bigl[ \frac1{\sqrt{\e}},\infty \bigr)$, when~$\alpha \leqslant d$, is 
\vspace{-.75\baselineskip}
\begin{align}
  \MoveEqLeft
\int_{\frac{1}{\sqrt{\e}}}^\infty  I_\e(r)^{\frac{p}2} r^{d-1}\,dr \leqslant \int_{\frac{1}{\sqrt{\e}}}^\infty  r^{d-1 - (d-\alpha)p} \biggl(\int_{\frac{\e r^2}{2(1+\e)}}^{\frac{r^2}2} t^{d-\alpha-1} e^{-t}\,dt\biggr)^{\frac{p}2}\,dr
  \\
 \label{e.lastpiece}
&\quad\leqslant  \frac1{(\alpha - d)^{\frac{p}2}}  \int_{\frac{1}{\sqrt{\e}}}^\infty  r^{d-1 - (d-\alpha)p} \exp\Bigl( - \frac{p r^2}{4}\Bigr) r^{(d-\alpha)p} \,dr 
= O\Bigl( e^{- \frac{p}{\e}} \Bigr)\,,
	\end{align}
which is, evidently, an exponentially small contribution. Inserting the computations from~\eqref{e.middlepiece},\eqref{e.firstpiece} and~\eqref{e.lastpiece} into~\eqref{e.computation-sub} we have shown that for~$\e$ small enough, 
\begin{equation}
	\|\varPhi_\e\|_{\mathscr{W}^{-\alpha,p}}^p \leqslant  2A_{d-1} \Gamma(d-\alpha)^{\frac{p}2}\begin{dcases}
	  \ln \paren[\Big]{\frac1{\e}} & \mbox{ if } \alpha p = d(p-1)\\
		\frac1{\alpha p - d(p-1)} & \mbox{ if } \alpha p > d(p-1)\,. 
	\end{dcases}
\end{equation}
As in Case 1, we must show a matching lower bound which demonstrates that these rates are asymptotically (for small~$\e$) optimal. To do this, we use~\eqref{e.decayiny} to compute in spherical coordinates, when~$\alpha p > d(p-1)$, that when~$\e$ is small enough,
\begin{align}
	\label{e.lowerboundintermediatealpha-1}
	&\mathscr{B}_0(\e) \geqslant A_{d-1}\int_{\sqrt{\e} \leqslant r \leqslant \frac1{\sqrt{\e}}} I_\e(r)^{\frac{p}2}r^{d-1} \,dr \\
&\quad 	\geqslant A_{d-1} \int_{\sqrt{\e}}^{\frac1{\sqrt{\e}}} \frac{r^{d-1}}{r^{p(d-\alpha)}} \biggl( \int_{ \frac{\e r^2}{2(1+\e)}}^{\frac{r^2}2} t^{d-\alpha-1} e^{-t}\,dt \biggr)^{\frac{p}2} \,dr \\
	&\quad = A_{d-1}\e^{\frac12 (\alpha p - (p-1)d)} \int_1^{\frac1{\e}} s^{\alpha p - (p-1)d -1} \biggl(  \int_{ \frac{\e^2 s^2}{2(1+\e)}}^{\frac{s^2}{2}}  t^{d-\alpha - 1} e^{-t}\,dt \biggr)^{\frac{p}2}\,ds \\
	&\quad \geqslant A_{d-1}\e^{\frac12 (\alpha p - (p-1)d)} \int_1^{\frac1{\e}} s^{\alpha p - (p-1)d -1} e^{- \frac{p s^2}{4}} \frac1{2(d- \alpha)^{p/2}} s^{(d-\alpha)p} \,ds  \\
	&\quad =  \frac{A_{d-1}}{2(d- \alpha)^{p/2}} \e^{\frac12 (\alpha p - (p-1)d)} \int_1^{\frac1{\e}} s^{d -1} e^{- \frac{p s^2}{4}} \,ds \\
	&\quad \geqslant   \frac{A_{d-1}}{4(d- \alpha)^{p/2}} \e^{\frac12 (\alpha p - (p-1)d)} \int_1^{\infty} s^{d -1} e^{- \frac{p s^2}{4}} \,ds\,. 
\end{align} 
In the previous calculation we made the change of variables~$r = \sqrt{\e}s$. 

In the critical case, that is, when~$\alpha p = d(p-1)$, we estimate slightly differently, arguing that, if~$d \geqslant p$, then
\begin{align}
	\label{e.lowerboundintermediatealpha-2}
	\mathscr{B}_0(\e) &\geqslant A_{d-1}\int_{\sqrt{\e} \leqslant r \leqslant \frac1{\sqrt{\e}}} I_\e(r)^{\frac{p}2} r^{d-1}\,dr \\
	& 	\geqslant A_{d-1} \int_{\sqrt{\e}}^{\frac1{\sqrt{\e}}} \frac{r^{d-1}}{r^{p(d-\alpha)}} \biggl( \int_{ \frac{\e r^2}{2(1+\e)}}^{\frac{r^2}2} t^{d-\alpha-1} e^{-t}\,dt \biggr)^{\frac{p}2} \,dr \\
	& = A_{d-1} \int_1^{\frac1{\e}} s^{ -1} \biggl(  \int_{ \frac{\e^2 s^2}{2(1+\e)}}^{\frac{s^2}{2}}  t^{\frac{d}p - 1} e^{-t}\,dt \biggr)^{\frac{p}2}\,ds \\
	& \geqslant A_{d-1} \int_1^{\frac1{\e}} s^{ -1} \biggl( \int_1^{\frac{s^2}2}  e^{-t}\,dt \biggr)^{\frac{p}2}\,ds\\
	& \geqslant A_{d-1}\int_1^{\frac1{\e}} s^{-1} \bigl(1 - e^{-\frac{s^2}2}\bigr)^{\frac{p}2} \,ds\\
	& \geqslant \bigl(1 - e^{-\frac12}\bigr)^{\frac{p}2} A_{d-1} \log \frac1{\e}\,, 
\end{align} 
whereas, when~$d < p$, 
\begin{align}
	\label{e.lowerboundintermediatealpha-3}
	\mathscr{B}_0(\e) &\geqslant A_{d-1}\int_{\sqrt{\e} \leqslant r \leqslant \frac1{\sqrt{\e}}} I_\e(r)^{\frac{p}2}r^{d-1} \,dr \\
	& 	\geqslant A_{d-1} \int_{\sqrt{\e}}^{\frac1{\sqrt{\e}}} \frac{r^{d-1}}{r^{p(d-\alpha)}} \biggl( \int_{ \frac{\e r^2}{2(1+\e)}}^{\frac{r^2}2} t^{d-\alpha-1} e^{-t}\,dt \biggr)^{\frac{p}2} \,dr \\
	& = A_{d-1} \int_1^{\frac1{\e}} s^{ -1} \biggl(  \int_{ \frac{\e^2 s^2}{2(1+\e)}}^{\frac{s^2}{2}}  t^{\frac{d}p - 1} e^{-t}\,dt \biggr)^{\frac{p}2}\,ds \\
	& \geqslant A_{d-1} \int_1^{\frac1{\e}} s^{ -1} \biggl( \int_{\frac{\e^2 s^2}{2(1+\e)}}^{1}  e^{-t}\,dt \biggr)^{\frac{p}2}\,ds\\
& \geqslant A_{d-1} \int_1^{\frac1{\e}} s^{ -1} \biggl( \int_{\frac12}^{1}  e^{-t}\,dt \biggr)^{\frac{p}2}\,ds\\
	& \geqslant \bigl(1 - e^{-\frac12}\bigr)^{\frac{p}2}A_{d-1}\int_1^{\frac1{\e}} s^{-1}  \,ds\\
	& \geqslant \bigl(1 - e^{-\frac12}\bigr)^{\frac{p}2} A_{d-1} \log \frac1{\e}\,. 
\end{align} 

Cases 1 and 2 together complete the proof of items (1) and (2) in Proposition~\ref{p:PhiEpBounds} for the quantity~$\mathscr{B}_0$. 

\smallskip 
  \case[$\frac{p-1}{p}d > \alpha$]
Finally, in this case, we compute 
 \begin{align} \label{e.computation-sub2}
 	\mathscr{B}_0(\e,\alpha,d,p)
	  &= 	A_{d-1} \int_0^\infty I_\e(r)^{\frac{p}2} r^{d-1}\,dr \\
 	& = A_{d-1}  \biggl( \int_0^{1} + \int_{1}^{\frac1{\sqrt{\e}}} + \int_{\frac1{\sqrt{\e}}}^\infty\biggr) I_\e(r)^{\frac{p}2} r^{d-1}\,dr \,. 
 \end{align}
 As in Case 1, we estimate these term-by-term. For the first piece we compute 
 \begin{equation*}
 	\label{e.firstpiecesubcrit}
 	\int_0^1 I_\e(r)^{\frac{p}2} r^{d-1}\,dr \leqslant \frac1{d} \Bigl( \frac1{\alpha} + \frac1{d-\alpha} \Bigr)^{\frac{p}2}\,,
 \end{equation*}
whereas for the middle piece we compute by plugging in~\eqref{e.gammabound} we obtain
\begin{equation} \label{e.middlepiece-subcrit}
	\int_{1}^{\frac1{\sqrt{\e}}} I_\e(r)^{\frac{p}2} r^{d-1} \,dr
	\quad \leqslant  \Gamma(d-\alpha)^{\frac{p}2} \int_{1}^{\frac1{\sqrt{\e}}}  \frac1{r^{p(d-\alpha)}} r^{d-1}\,dr	\leqslant  \frac{\Gamma(d-\alpha)^{\frac{p}2}}{d(p-1) - \alpha p}\,. 
	\end{equation}
The estimate for the last piece is exactly the same as~\eqref{e.lastpiece}. To conclude, by the monotone convergence theorem we obtain that in this sub-critical case, 
\begin{equation*}
	\lim_{\e \to 0^+} \mathscr{B}_0(\e,\alpha,d,p) =\int_{\Rd} \biggl(\int_0^{\infty} \frac{s^{\alpha-1}}{ (1+s)^d} \exp \biggl(- \frac{|y|^2}{2(1+s)}\biggr)\,ds\biggr)^{\frac{p}2}\,dy \in (0,\infty)\,. 
	\qedhere
\end{equation*}
\end{proof}

\bibliographystyle{halpha-abbrv}
\bibliography{gautam-refs1,gautam-refs2,raghav-refs,preprints}
\end{document}